\documentclass[reqno]{amsart}
\usepackage{amssymb, amsmath, amsthm, amscd, mathrsfs}

\usepackage{paralist}
 
\usepackage{bigints}
\usepackage{dsfont}
\usepackage{empheq}
\usepackage{amssymb}
\usepackage{cases}
\usepackage{enumitem}
\usepackage{cite}
\allowdisplaybreaks

\usepackage[colorlinks]{hyperref}
\usepackage[colorinlistoftodos]{todonotes}

\theoremstyle{plain}
\newtheorem{theorem}{Theorem}[section]

\newtheorem{lemma}[theorem]{Lemma}
\newtheorem{proposition}{Proposition}

\theoremstyle{definition}

\newtheorem{remark}{Remark}

\def \dv {\mathrm{div}}
\def \d {\mathrm{d}}

\title[Boundary null controllability for the heat equation] 
      {Boundary null controllability for the heat equation with dynamic boundary conditions}

\author{S. E. Chorfi}
\author{G. El Guermai}
\author{A. Khoutaibi}
\author{L. Maniar}

\address{S. E. Chorfi, G. El Guermai, A. Khoutaibi and L. Maniar, Cadi Ayyad University, Faculty of Sciences Semlalia, LMDP, UMMISCO (IRD-UPMC), B.P. 2390, Marrakesh, Morocco}

\email{chorphi@gmail.com, ghita.el.guermai@gmail.com, akhoutaibi2@gmail.com, maniar@uca.ma}

\subjclass[2010]{Primary: 93B05; Secondary: 35K05, 93B07.}
 \keywords{Null controllability, parabolic equations, Carleman estimate, dynamic boundary conditions, surface diffusion.}

\begin{document}
\begin{abstract}
In this paper, we are concerned with the boundary controllability of heat equation with dynamic boundary conditions. More precisely, we prove that the equation is null controllable at any positive time by means of a boundary control supported on an arbitrary subboundary. The proof of the main result combines a new boundary Carleman estimate and some regularity estimates for the adjoint system, with an explicit dependence with respect to the final time. This technique allows us to overcome a new difficulty that arises when absorbing a normal derivative term.
\end{abstract}

\maketitle

\section{Introduction}

Let us fix $T>0$ and an integer $N \ge 2$. Throughout the paper, $\Omega \subset \mathbb{R}^N$ is a bounded domain with boundary $\Gamma=\partial \Omega$ of class $C^{3,1}$ and $\Gamma^0 \subset \Gamma$ a nonempty relatively open subset of $\Gamma$. Let us set
$$E_{T}=(0, T) \times E \qquad \text{ for all } E \subset \overline{\Omega}.$$
Then, we consider the following heat equation with dynamic boundary conditions
\begin{empheq}[left = \empheqlbrace]{alignat=2} \label{eq:intr}
\begin{aligned}
&\partial_t y - d\Delta y = f(t,x), &\qquad\text{ in } \Omega_T , \\
&\partial_t y_{\Gamma} -\delta\Delta_\Gamma y_\Gamma + d\partial_\nu y = g(t,x), &\qquad\text{ on } \Gamma_T , \\
& y_{\Gamma}(t,x) = y_{|\Gamma}(t,x), &\qquad\text{ on } \Gamma_T , \\
& (y,y_\Gamma)\rvert_{t=0}=(y_0,y_{0,\Gamma}), &\Omega\times\Gamma,
\end{aligned}
\end{empheq}
where $d, \delta >0$ are positive real numbers, $f \in L^{2}(\Omega_T)$ and $g \in L^{2}(\Gamma_T)$ are the source terms, $\partial_{\nu} y:=(\nabla y \cdot \nu)_{\mid \Gamma}$ is the normal derivative, where $\nu$ is the outer unit normal field on $\Gamma$ exterior to $\Omega$. Further, we denote by $y_{\Gamma}:=y_{\mid \Gamma}$ the trace of $y$. Finally,
$\Delta$ is the Laplace operator, $\Delta_{\Gamma}$ is the Laplace-Beltrami operator on the Riemannian submanifold $\Gamma$ (see \cite{Ta'11}), and $y_{0} \in L^{2}(\Omega), y_{0,\Gamma} \in L^{2}(\Gamma)$ are initial states, which are not necessarily interrelated. The term $d \partial_{\nu} y$ represents the interaction domain-boundary, while $\delta \Delta_{\Gamma} y_{\Gamma}$ stands for a boundary diffusion. 

The aim of this paper is to prove the boundary null controllability of the system \eqref{eq:intr}. Our main result relies on a new boundary Carleman estimate for the adjoint system associated to \eqref{eq:intr}, with a local integral in terms of Dirichlet data localized in the arbitrary subboundary $\Gamma^0 \subset \Gamma$.

Evolution equations with dynamic boundary conditions have attracted a considerable attention in last years, and have found applications in various mathematical models such as population dynamics and heat transfer, in particular when the diffusion process takes place in the interface between a solid and a moving fluid, see \cite{TMR'14, FH'11, La'32} and the references therein. The reader can be referred to the seminal paper \cite{Go'06} for the physical interpretation and derivation of such boundary conditions. We also refer to the recent paper \cite{Sa'20}, where the derivation is given by employing the Carslaw-Jaeger relation. Many researches were done for the wellposedness and regularity results related to such equations, see \cite{FGGR'02, Go'06, MMS'17, VV'11} for the case of a scalar diffusion and without gradient terms, and recently \cite{KMMR'19}, in the presence of gradient terms. The general case of variable diffusion coefficients with gradient terms has been recently studied in \cite{ACMO'20}. Although the controllability of parabolic equations, by the Carleman estimate approach, is well studied in the case of static boundary conditions (Dirichlet, Neumann or Fourier boundary conditions), see, e.g., \cite{CIK'96, FG'06, FGGP'06, FI'96, Im'95, IY'01, IY'98}, only some recent works concern the controllability and inverse problems in the dynamic boundary conditions case; see \cite{ACM'21, ACM'21', BCMO'20}.

In \cite{KM'19, MMS'17}, the authors have considered the distributed null controllability for \eqref{eq:intr}, that is, the null controllability from an arbitrary subdomain $\omega \Subset \Omega$. The inverse source problem with interior measurements has been recently studied in a more general setting in \cite{ACMO'20}. In these works, we have adopted the interior weight function used in \cite{FI'96} for parabolic equations with classical boundary conditions. The aforementioned weight function is null on the boundary of the domain. This fact along with the nature of classical boundary conditions allow several boundary terms to be canceled. This is not the case in the present setting.

As for the standard (scalar) heat equation, a well-known result states that the boundary null controllability and distributed null controllability are equivalent (see \cite[Theorem 2.2]{ABBT'11}). This result is no longer valid for coupled parabolic systems. In our case, one can obtain the boundary null controllability from a distributed control by enlarging the domain $\Omega$. Nevertheless, this technique for dynamic boundary conditions requires more regular initial data than those in the Dirichlet or Neumann cases, see \cite[Theorem 4.5]{MMS'17}. More precisely, due to the use of Sobolev's embedding and local regularity results, the authors take initial conditions such that $y_0 \in H^2(\Omega), \, {y_0}_{|\Gamma}\in H^2(\Gamma)$ and $y_0\in W^{2-\frac{2}{p}}_p(\Omega)$ for some $p>1+\frac{N}{2}$. Our boundary Carleman estimate allows us to recover the boundary controllability result without assuming any further regularity assumption. Furthermore, our direct approach has the following extra advantages:
\begin{itemize}
\item[$\bullet$] we keep an explicit dependence of the constants in the boundary Carleman estimate with respect to the final time $T$. This fact is crucial when dealing with semilinear problems;
\item[$\bullet$] we obtain an explicit constant, with respect to the final time $T$, of the boundary observability constant which is different from the interior observability constant obtained in \cite{KM'19};
\item[$\bullet$] we then obtain a continuous dependence of null controls and corresponding states with respect to initial data and source terms;
\item[$\bullet$] our boundary Carleman estimate allows to obtain an explicit upper bound for the cost of boundary approximate controllability, i.e., an explicit bound for the size of the minimal norm control. We refer to \cite{BCMO'20} for the case of a distributed control;
\item[$\bullet$] the techniques could seemingly be extended to deal with inverse problems with boundary measurements in the spirit of \cite{ACM'21, ACMO'20}. However, a careful analysis is needed to include a large class of source terms.
\end{itemize}

In comparison with the classical parabolic equations, where the control acts by means of the Dirichlet condition on a part of the boundary of the domain, the present Carleman inequality requires a different technique to absorb a term of Neumann type (in terms of the normal derivative) by the left-hand side of the Carleman estimate. Absorbing this term turns out to be a main difficulty that will be overcome with help of some parabolic regularity estimates with explicit constants in terms of the final time $T$ (see Proposition \ref{prop:accreg}). Such regularity estimates can be very useful while doing Carleman estimates \cite{BGP'19}. It is worth mentioning the paper \cite{CG'12} which is the first to deal with a Carleman inequality when a time derivative appears in the boundary condition (but without a boundary diffusion). However, it only consider the one dimensional case and point out the difficulty of the controllability in higher dimensions for such equations (see Section 4 in \cite{CG'12}).

In the subsequent sections, we typically use the following surface divergence theorem
\begin{equation} \label{eq:sdt}
 \int_{\Gamma} \Delta_\Gamma y\,z \, \d S 
= - \int_{\Gamma } \langle\nabla_\Gamma y, \nabla_\Gamma z\rangle_\Gamma \, \d S, \qquad 
      y \in H^2(\Gamma), \;z\in H^1(\Gamma), 
\end{equation}
where $\nabla_\Gamma$ is the tangential gradient and $\langle\cdot, \cdot\rangle_\Gamma$ is the Riemannian inner product of tangential vectors on $\Gamma$ \cite{Jo'08}. To simplify notation, we will denote $$\left \langle \nabla_{\Gamma}f, \nabla_{\Gamma}g\right \rangle_{\Gamma} =  \nabla_{\Gamma}f\cdot  \nabla_{\Gamma}g.$$
Also, we need the following identity which holds for general Riemannian manifolds. Its proof can be found in \cite{BY'17}:
\begin{align} \label{formula:grad}
\nabla_{\Gamma}(\nabla_{\Gamma}f\cdot \nabla_{\Gamma}g)\cdot\nabla_{\Gamma}f=  \nabla_{\Gamma}^{2} g\left(\nabla_{\Gamma} f, \nabla_{\Gamma} f\right)+\frac{1}{2} \nabla_{\Gamma} g \cdot \nabla_{\Gamma}\left(\left|\nabla_{\Gamma} f\right|^{2}\right)
\end{align}
for $f\in H^2(\Gamma)$ and $g\in H^2(\Gamma)$, where $\nabla_{\Gamma}^{2}$ denotes the Hessian with respect to the Riemannian metric on $\Gamma$.

The remainder of the paper is organized as follows. In Section \ref{sec:well}, we recall some wellposedness and regularity results that will be useful in the sequel. In Section \ref{sec:carleman}, we prove a boundary Carleman estimate for the homogeneous adjoint system. Finally, in Section \ref{sec:null}, we apply the Carleman estimate to show the observability inequality, and consequently we obtain the boundary null controllability of the system \eqref{eq:intr}.

\section{Wellposedness and regularity of the solution} \label{sec:well}
In this section, we recall some wellposedness and regularity results related to system \eqref{eq:intr}. The reader can refer to \cite{MMS'17} for different notions of the solution as well as detailed proofs.

\subsection{Notations}
First, we introduce the following real Hilbert space
$$\mathbb{L}^2:=L^2(\Omega, \d x)\times L^2(\Gamma, \d S),$$
equipped with the scalar product given by
$$\langle (y,y_\Gamma),(z,z_\Gamma)\rangle_{\mathbb{L}^2} =\langle y,z\rangle_{L^2(\Omega)} +\langle y_\Gamma,z_\Gamma\rangle_{L^2(\Gamma)},$$
where we denoted the Lebesgue measure on $\Omega$ and the surface measure on $\Gamma$ by $\d x$ and $\d S$, respectively. We also define the spaces
$$\mathbb{H}^k:=\left\{(y,y_\Gamma)\in H^k(\Omega)\times H^k(\Gamma)\colon y_{|\Gamma} =y_\Gamma \right\}, \quad k=1,2, \ldots,$$ which are Hilbert spaces when equipped with the standard product norm. Note that $\mathbb{H}^1$ is a Hilbert space when equipped with the norm
\begin{align*}
\left\|\left(y, y_{\Gamma}\right)\right\|_{\mathbb{H}^{1}}&=\left\langle\left(y, y_{\Gamma}\right),\left(y, y_{\Gamma}\right)\right\rangle_{\mathbb{H}^{1}}^{1 / 2},\\
\text{where }\left\langle\left(y, y_{\Gamma}\right),\left(z, z_{\Gamma}\right)\right\rangle_{\mathbb{H}^{1}}&=\int_{\Omega} \nabla y \cdot\nabla z \,\d x +\int_{\Gamma} \nabla_{\Gamma} y_{\Gamma} \cdot\nabla_{\Gamma} z_{\Gamma} \,\d S +\int_{\Gamma} y_{\Gamma} z_{\Gamma} \,\d S.
\end{align*}
Similarly, $\mathbb{H}^2$ is a Hilbert space with respect to the norm
\begin{align*}
\left\|\left(y, y_{\Gamma}\right)\right\|_{\mathbb{H}^{2}}&=\left\langle\left(y, y_{\Gamma}\right),\left(y, y_{\Gamma}\right)\right\rangle_{\mathbb{H}^{2}}^{1 / 2},\\
\text{where }\left\langle\left(y, y_{\Gamma}\right),\left(z, z_{\Gamma}\right)\right\rangle_{\mathbb{H}^{2}}&=\int_{\Omega} \Delta y \Delta z \,\d x +\int_{\Gamma} \Delta_{\Gamma} y_{\Gamma} \Delta_{\Gamma} z_{\Gamma} \,\d S +\int_{\Gamma} y_{\Gamma} z_{\Gamma} \,\d S.
\end{align*}
For the regularity of the solutions, we introduce the following spaces
$$\mathbb{E}_1(t_0,t_1):=H^1\left(t_0,t_1 ;\mathbb{L}^2\right) \cap L^2\left(t_0,t_1 ;\mathbb{H}^2\right) \text{  for } t_1 >t_0 \text{ in } \mathbb{R},$$
$$\mathbb{E}_2(t_0,t_1):=H^1\left(t_0,t_1;\mathbb{H}^2\right) \cap H^2\left(t_0,t_1;\mathbb{L}^2\right) \text{  for } t_1 >t_0 \text{ in } \mathbb{R}.$$
In particular,
$$\mathbb{E}_1 := \mathbb{E}_1(0,T) \;\text{ and }\; \mathbb{E}_2:= \mathbb{E}_2(0,T).$$

Finally, we recall that the system \eqref{eq:intr} is governed by the linear operator that is denoted by
$$\mathcal{A} =\left ( \begin{array}{cc} d\Delta & 0 \\ -d\partial_\nu & \delta \Delta_\Gamma
\end{array}\right), \qquad D(\mathcal{A}) = \mathbb{H}^2.
$$
It is shown in \cite[Proposition 2.1]{MMS'17} that the operator $\mathcal{A}$ is self-adjoint dissipative and generates an analytic $C_0$-semigroup $\left(\mathrm{e}^{t\mathcal{A}}\right)_{t\geq 0}$ on $\mathbb{L}^2$.

\subsection{Regularity estimates}
We first recall an elliptic regularity result, which will allow us to characterize the domain of $\mathcal{A}^2$ and prove parabolic regularity later on. Consider the following elliptic problem
\begin{equation}\label{ell:prob}
\begin{cases}
d\Delta u =f, &\mbox{ in }\Omega,\\
\delta \Delta_\Gamma u_\Gamma - d\partial_\nu u - u_\Gamma =g, & \mbox{ on }\Gamma, \\
u_\Gamma=u_{|\Gamma}, & \mbox{ on }\Gamma.
\end{cases}
\end{equation}
The following result is a special case of \cite[Theorem 3.4]{KCDQ'15}.
\begin{proposition}
If $\mathcal{F}=(f,g)\in H^2(\Omega)\times H^2(\Gamma)$, then the solution of \eqref{ell:prob} satisfies $(u,u_\Gamma)\in \mathbb{H}^4$ and
\begin{equation}\label{eq:norm}
\|(u,u_\Gamma)\|_{\mathbb{H}^{4}} \leq C \|\mathcal{F}\|_{H^2(\Omega)\times H^2(\Gamma)},
\end{equation}
where $C=C(d,\delta)$ is a positive constant.
\end{proposition}
Consequently, we have the following lemma.
\begin{lemma}\label{lem:dom}
The domain of $\mathcal{A}^2$ is given by
$$D\left(\mathcal{A}^2\right)=\left\{(u,u_\Gamma) \in \mathbb{H}^4 \colon (d\Delta u)\rvert_\Gamma =\delta \Delta_\Gamma u_\Gamma - d\partial_\nu u\right\},$$
which is a Hilbert space equipped with the standard inner product of $H^4(\Omega)\times H^4(\Gamma)$.
\end{lemma}
Now, let us consider the backward adjoint problem of \eqref{eq:intr}, which is given by
\begin{empheq}[left = \empheqlbrace]{alignat=2}
\begin{aligned} \label{eq:adj}
&-\partial_t\varphi-d\Delta \varphi = 0, &\qquad\text{ in } \Omega_T , \\
&-\partial_t \varphi_\Gamma -\delta\Delta_\Gamma \varphi_\Gamma +d\partial_\nu \varphi  = 0, &\qquad\text{ on } \Gamma_T , \\
& \varphi_{\Gamma}(t,x) = \varphi_{|\Gamma}(t,x), &\qquad\text{ on } \Gamma_T , \\
& (\varphi(T,\cdot),\varphi_\Gamma(T,\cdot)) = (\varphi_T, \varphi_{T,\Gamma}), &\Omega\times\Gamma
\end{aligned}
\end{empheq}
for given $(\varphi_T, \varphi_{T,\Gamma})$ in $\mathbb{L}^2$. We recall the following result from \cite{MMS'17} for future use.
\begin{proposition}\label{prop:mild}
Let $f\in L^2(\Omega_T)$, $g\in L^2(\Gamma_T)$, $Y_0:=(y_0,y_{0,\Gamma})\in \mathbb{L}^2$ and $(\varphi_T,\varphi_{T,\Gamma})\in \mathbb{L}^2$. Then the following assertions are true.
\begin{itemize}
 \item[\text{(a)}] There is a unique mild solution $Y:=(y,y_\Gamma)\in C([0,T];\mathbb{L}^2)$ of \eqref{eq:intr}. Moreover, $Y$ belongs to $\mathbb{E}_1(\tau,T)$ and solves \eqref{eq:intr} strongly on $(\tau,T)$ with initial data $Y(\tau)$, for each $\tau\in(0,T)$. This mild solution is given by
\begin{equation*}
Y(t)= \mathrm{e}^{t\mathcal{A}} Y_0 +\int_0^t \mathrm{e}^{(t-s)\mathcal{A}}(f(s),g(s))\, \d s, \qquad t\in [0,T].
\end{equation*}
\item[\text{(b)}] There is a constant $C>0$ such that for all data the mild solution $Y$ satisfies
\begin{equation} \label{eq:esti-1}
 \|Y\|_{C([0,T];\mathbb{L}^2)} \leq C \left(\|Y_0\|_{\mathbb{L}^2} + \|f\|_{L^2(\Omega_T)} + \|g\|_{L^2(\Gamma_T)}\right). 
\end{equation}
 \item[\text{(c)}] If $Y_0\in \mathbb{H}^1$, then the mild solution from (a) is the strong one.
\item[\text{(d)}] The backward adjoint system \eqref{eq:adj} admits a unique mild solution given by $(\varphi(t),\varphi_\Gamma(t))=\mathrm{e}^{(T-t)\mathcal{A}}(\varphi_T,\varphi_{T,\Gamma})$. It is the unique strong solution if $(\varphi_T,\varphi_{T,\Gamma})\in \mathbb{H}^1$. 
\end{itemize}
\end{proposition}

In order to absorb the normal derivative term in the proof of the Carleman inequality and obtain an explicit observability constant in terms of $T$, we shall prove more accurate regularity estimates for the following adjoint system
\begin{equation}\label{acp1}
\left\{\begin{array}{l}
\hspace{-0.19cm} -\partial_t \Phi(t)=\mathcal{A} \Phi(t) + \mathcal{F}(t), \quad t \in (0, T], \\
\hspace{-0.19cm} \Phi(T)=(0,0),
\end{array}\right.
\end{equation}
where $\Phi:=(\varphi,\varphi_\Gamma)$ and $\mathcal{F}:=(f,g)$.

\begin{proposition}\label{prop:accreg}
The following assertions hold
\begin{itemize}
\item[(i)] if $\mathcal{F}\in L^2\left(0,T;\mathbb{L}^2\right)$, then the solution $\Phi$ of \eqref{acp1} satisfies 
\begin{equation} \label{I0}
\begin{aligned}
& \quad\qquad \Phi \in L^2\left(0,T;\mathbb{H}^1\right),\\
&\|\Phi\|_{L^2\left(0,T;\mathbb{H}^1\right)} \le C (1+T) \|\mathcal{F}\|_{L^2\left(0,T;\mathbb{L}^2\right)}, 
\end{aligned}
\end{equation}
where $C>0$ is independent of $T$ and $\mathcal{F}$.
\item[(ii)] If $\mathcal{F}\in L^2\left(0,T;\mathbb{L}^2\right)$, then the solution $\Phi$ of \eqref{acp1} satisfies 
\begin{equation}\label{I2}
\begin{aligned}
& \quad\qquad \Phi \in L^2\left(0,T;\mathbb{H}^2\right)\cap H^1\left(0,T; \mathbb{L}^2\right),\\
&\|\Phi\|_{L^2\left(0,T;\mathbb{H}^2\right)} + \|\partial_t \Phi\|_{L^2\left(0,T;\mathbb{L}^2\right)} \le C (1+T) \|\mathcal{F}\|_{L^2\left(0,T;\mathbb{L}^2\right)}, 
\end{aligned}
\end{equation}
\item[(iii)] if $\mathcal{F}\in L^2\left(0,T;\mathbb{H}^2\right)\cap H^1\left(0,T; \mathbb{L}^2\right)$ and $\mathcal{F}(T)=(0,0)$, then the solution $\Phi$ of \eqref{acp1} satisfies 
\begin{equation}
\begin{aligned}
& \hspace{2.5cm} \Phi \in L^2\left(0,T;\mathbb{H}^4\right)\cap H^2\left(0,T; \mathbb{L}^2\right),\\
&\|\Phi\|_{L^2\left(0,T;\mathbb{H}^4\right)} + \left\|\partial_t^2 \Phi\right\|_{L^2\left(0,T;\mathbb{L}^2\right)} \le C (1+T) \left(\|\mathcal{F}\|_{L^2\left(0,T;\mathbb{H}^2\right)} + \|\partial_t \mathcal{F}\|_{L^2\left(0,T; \mathbb{L}^2\right)} \right).
\end{aligned}
\end{equation}
\end{itemize}
\end{proposition}

\begin{proof}
By density, it suffices to consider a smooth source term so that the solution $\Phi$ is smooth. Next, we focus on the regularity estimates.

$(i)$ Multiplying \eqref{eq:adj}$_{1}$ by $\varphi$ and \eqref{eq:adj}$_{2}$ by $\varphi_{\Gamma}$ and using Green formula and surface diverge formula, we obtain, for any $t\in (0,T)$, that 
\begin{equation}\label{eg:3}
    -\frac{1}{2}\partial_{t}\|\Phi (t)\|_{\mathbb{L}^2}^2+d\int_{\Omega}|\nabla\varphi(t)|^2 \,\d x+\delta\int_{\Gamma}|\nabla_{\Gamma}\varphi_{\Gamma}(t)|^2 \, \d S=\left\langle \mathcal{F}(t), \Phi(t)\right\rangle _{\mathbb{L}^2}.
\end{equation}
Integrating over $(t,T)$, we infer, for any $\varepsilon>0$ and any $t\in (0,T)$, that
\begin{align*}
\|\Phi(t)\|_{\mathbb{L}^2}^2 &\le 2 \int_0^T \left\langle \mathcal{F}(t), \Phi(t)\right\rangle _{\mathbb{L}^2}\,\d t \\
& \le 2 \int_0^T \|\mathcal{F}(t)\|_{\mathbb{L}^2} \|\Phi(t)\|_{\mathbb{L}^2}\,\d t \\
& \le 2 \|\mathcal{F}\|_{L^1(0,T;\mathbb{L}^2)} \|\Phi\|_{C\left([0,T];\mathbb{L}^2\right)}\\
& \le C_\varepsilon \|\mathcal{F}\|_{L^1(0,T;\mathbb{L}^2)}^2 + \varepsilon \|\Phi\|_{C([0,T];\mathbb{L}^2)}^2.
\end{align*}
Choosing $\varepsilon>0$ small enough, we find that
$$\|\Phi\|_{C\left([0,T];\mathbb{L}^2\right)}^2 \le C \|\mathcal{F}\|_{L^1(0,T;\mathbb{L}^2)}^2.$$
Since
\begin{equation*}
\|\Phi\|_{L^2(0,T;\mathbb{L}^2)}^2 \le T \|\Phi\|_{C\left([0,T];\mathbb{L}^2\right)}^2 \qquad \text{ and }\qquad \|\mathcal{F}\|_{L^1(0,T;\mathbb{L}^2)}^2\leq T\|\mathcal{F}\|_{L^2(0,T;\mathbb{L}^2)}^2,
\end{equation*}
it follows that
\begin{equation}\label{eg:6}
\|\Phi\|_{L^2(0,T;\mathbb{L}^2)} \le C T \|\mathcal{F}\|_{L^2(0,T;\mathbb{L}^2)}.
\end{equation}
On the other hand, integrating \eqref{eg:3} over $(0,T)$, we obtain
\begin{equation}\label{eg:5}
\min(d,\delta) \|(\nabla\varphi, \nabla_{\Gamma}\varphi_{\Gamma})\|_{L^2(0,T;\mathbb{L}^2)}^2 \le \int_{0}^T\left\langle \mathcal{F}(t),\Phi(t)\right\rangle_{\mathbb{L}^2}\, \d t.
\end{equation}
We have 
\begin{equation*}
    2\int_{0}^T\left\langle \mathcal{F}(t),\Phi(t)\right\rangle_{\mathbb{L}^2} \,\d t\leq \|\Phi\|_{L^2(0,T;\mathbb{L}^2)}^2 + \|\mathcal{F}\|_{L^2(0,T;\mathbb{L}^2)}^2.
\end{equation*}
This inequality with \eqref{eg:6}-\eqref{eg:5}, yield
$$
\|\Phi\|_{L^2(0,T;\mathbb{H}^1)}\leq C(1+T)\|\mathcal{F}\|_{L^2(0,T;\mathbb{L}^2)}.
$$
\smallskip

$(ii)$ For $t\in (0,T)$, we have
\begin{align*}
\int_\Omega f^2 \,\d x &=\int_\Omega (\partial_t \varphi +d\Delta \varphi)^2 \,\d x \notag\\
& = \int_\Omega |\partial_t \varphi|^2 \,\d x - d \int_\Omega \partial_t |\nabla \varphi|^2 \,\d x + 2d \int_\Gamma \partial_t \varphi_\Gamma \partial_\nu \varphi \,\d S + d^2 \int_\Omega |\Delta \varphi|^2 \,\d x, \\
\int_\Gamma g^2 \,\d S &=\int_\Gamma (\partial_t \varphi_\Gamma + \delta\Delta_\Gamma \varphi_\Gamma - d\partial_\nu \varphi)^2 \,\d S \notag\\
& \hspace{-0.6cm} = \int_\Gamma |\partial_t \varphi_\Gamma|^2 \,\d S - \delta \int_\Gamma \partial_t |\nabla_\Gamma \varphi_\Gamma|^2 \,\d S - 2d \int_\Gamma \partial_t \varphi_\Gamma \partial_\nu \varphi \,\d S + \delta^2 \int_\Gamma |\Delta_\Gamma \varphi_\Gamma|^2 \,\d S \notag\\
& - 2 d\delta \int_\Gamma \Delta_\Gamma \varphi_\Gamma \partial_\nu \varphi\, \d S,
\end{align*}
where we used Green's formula and surface divergence formula. Adding up the two previous identities and integrating over $(0,T)$, we obtain
\begin{align*}
\|\Phi\|_{L^2\left(0,T;\mathbb{H}^2\right)}^2 + \|\partial_t \Phi\|_{L^2\left(0,T;\mathbb{L}^2\right)}^2 &\le C(d,\delta) \left(\|\mathcal{F}\|_{L^2\left(0,T;\mathbb{L}^2\right)}^2 + \|\varphi_\Gamma\|_{L^2(0,T;L^2(\Gamma))}^2\right)\\
& \quad + 2 d\delta \int_0^T \int_\Gamma |\Delta_\Gamma \varphi_\Gamma| |\partial_\nu \varphi| \,\d S.
\end{align*}
The second term on the right-hand side can be estimated by \eqref{I0}. For the last term, we use Young's inequality for $\varepsilon>0$ to write
$$2 \int_0^T \int_\Gamma |\Delta_\Gamma \varphi_\Gamma| |\partial_\nu \varphi| \,\d S \le \varepsilon \int_0^T \int_\Gamma |\Delta_\Gamma \varphi_\Gamma|^2 \, \d S + \frac{1}{\varepsilon} \|\partial_\nu \varphi\|_{L^2(0,T;L^2(\Gamma))}^2.$$
Fixing $\varepsilon$ small enough, the first term on the right-hand side can be estimated by the left-hand side of the penultimate estimate.

Next, it suffices to estimate the normal derivative term so that
\begin{equation*}
      \|\partial_{\nu}\varphi\|_{L^2(\Gamma)}^2\leq \varepsilon\|\partial_{t}\varphi \|_{L^2(\Omega)}^2+\|\nabla_{\Gamma}\varphi\|_{L^2(\Gamma)}^2 +C_\varepsilon\|\nabla \varphi\|_{L^2(\Omega)}^2+\|f\|_{L^2(\Omega)}^2,
\end{equation*}
where $\varepsilon>0$ and C is a positive constant independent of T. Indeed, since $\Omega$ is of class $C^2$, there exists $\mathbf{N}\in C^2\left(\overline{\Omega}\right)^N$ (see \cite[Lemma 2.3]{BY'17}) such that
\[ 
\mathbf{N}(x)=\nu(x),\; x\in \Gamma, \qquad \text{ and }\qquad |\mathbf{N}(x)|\leq 1, \; x\in \Omega.
\]
Multiplying \eqref{eq:adj}$_{1}$ by $\nabla \varphi \cdot \mathbf{N}$ and integrating over $\Omega$, we obtain 
\begin{equation}\label{eg:1}
    -\int_{\Omega}\partial_{t}\varphi \nabla \varphi\cdot \mathbf{N} \,\d x-\int_{\Omega}\Delta\varphi \nabla\varphi \cdot \mathbf{N} \,\d x=\int_{\Omega} f \nabla \varphi\cdot \mathbf{N} \,\d x.
\end{equation}
By Green's formula, we obtain
\begin{align}\label{eg:2}
    -\int_{\Omega}\Delta \varphi \nabla \varphi\cdot \mathbf{N} \,\d x &=-\int_{\Gamma}|\partial_{\nu}\varphi|^2 \,\d S+\int_{\Omega}\nabla \varphi\cdot\nabla(\nabla\varphi\cdot \mathbf{N}) \,\d x\nonumber\\
    &=-\int_{\Gamma}|\partial_{\nu}\varphi|^2 \,\d S+\frac{1}{2}\int_{\Omega}\nabla |\nabla\varphi|^2 \cdot \mathbf{N} \,\d x +\int_{\Omega} \nabla \mathbf{N} \nabla \varphi \cdot \nabla \varphi \,\d x\nonumber\\
    &=-\int_{\Gamma}|\partial_{\nu}\varphi|^2\,\d S+\frac{1}{2}\int_{\Gamma}|\nabla\varphi|^2\,\d S-\frac{1}{2}\int_{\Omega} \dv\,\mathbf{N}\;|\nabla\varphi|^2\,\d x\nonumber\\
    & \qquad +\int_{\Omega} \nabla \mathbf{N} \nabla \varphi \cdot \nabla \varphi \,\d x \nonumber\\
    &=-\frac{1}{2}\int_{\Gamma}|\partial_{\nu}\varphi|^2\,\d S+\frac{1}{2}\int_{\Gamma}|\nabla_{\Gamma}\varphi|^2\,\d S-\frac{1}{2}\int_{\Omega} \dv\,\mathbf{N}\; |\nabla\varphi|^2\,\d x\nonumber\\
    &\qquad +\int_{\Omega} \nabla \mathbf{N} \nabla \varphi \cdot \nabla \varphi \,\d x,
\end{align}
where $\nabla \mathbf{N} =\left(\partial_{i}\mathbf{N}_{j}\right)$ stands for the standard Jacobian matrix of $\mathbf{N}$.

Using \eqref{eg:1} and \eqref{eg:2}, we obtain
\begin{align*}
\|\partial_{\nu}\varphi\|_{L^2(\Gamma)}^2&=-2\int_{\Omega}\partial_{t}\varphi\nabla\varphi\cdot \mathbf{N}\, \d x+\|\nabla_{\Gamma}\varphi\|_{L^2(\Gamma)}^2-\int_{\Omega} \dv\, \mathbf{N}\; |\nabla\varphi|^2\, \d x\\
& \qquad +2\int_{\Omega} \nabla \mathbf{N} \nabla \varphi \cdot \nabla \varphi \,\d x-2\int_{\Omega} f \nabla\varphi\cdot \mathbf{N} \, \d x.
\end{align*}
Let $\varepsilon>0$. Using Young's inequality and the fact $\mathbf{N}\in C^2\left(\overline{\Omega}\right)^N$, there exists a positive constant $C$ independent of $T$ such that
\begin{equation*}
      \|\partial_{\nu}\varphi\|_{L^2(\Gamma)}^2\leq \varepsilon\|\partial_{t}\varphi \|_{L^2(\Omega)}^2+\|\nabla_{\Gamma}\varphi\|_{L^2(\Gamma)}^2 +C_\varepsilon\|\nabla \varphi\|_{L^2(\Omega)}^2+\|f\|_{L^2(\Omega)}^2.
\end{equation*}
This completes the proof of $(ii)$.
\smallskip

$(iii)$ We set $\Psi=\partial_t \Phi$. Differentiating \eqref{acp1} in time $t$, we see that $\Psi$ solves the same system \eqref{acp1} with source term $\partial_t\mathcal{F}$. Then by (ii), we have
\begin{equation}\label{I3}
\|\partial_t \Phi\|_{L^2\left(0,T;\mathbb{H}^2\right)} + \left\|\partial_t^2 \Phi\right\|_{L^2\left(0,T;\mathbb{L}^2\right)} \le C(1+T) \|\partial_t\mathcal{F}\|_{L^2\left(0,T;\mathbb{L}^2\right)}.
\end{equation}
On the other hand, $(\varphi, \varphi_\Gamma)$ solves the following problem
\begin{equation*}
\begin{cases}
d\Delta \varphi =-f -\partial_t\varphi, &\mbox{ in }\Omega,\\
\delta \Delta_\Gamma \varphi_\Gamma - d\partial_\nu \varphi - \varphi_\Gamma =-g-\partial_t\varphi_\Gamma -\varphi_\Gamma , & \mbox{ on }\Gamma, \\
\varphi_\Gamma=\varphi_{|\Gamma}, & \mbox{ on }\Gamma.
\end{cases}
\end{equation*}
Then inequality \eqref{eq:norm} implies that
\begin{equation*}
\|\Phi\|_{\mathbb{H}^{4}}^2 \leq C \left( \|\mathcal{F}\|_{\mathbb{H}^2}^2 + \|\partial_t\Phi\|_{\mathbb{H}^2}^2 + \|\Phi\|_{\mathbb{H}^2}^2\right).
\end{equation*}
Integrating over $(0,T)$ and combining the resulting inequality with \eqref{I2} and \eqref{I3}, we obtain the desired estimate.
\end{proof}

\begin{remark}
We emphasize that a similar regularity result to Proposition \ref{prop:accreg} can be obtained by classical semigroup theory and Lemma \ref{lem:dom}, see, for instance, Theorems 2.3 and 2.4 in \cite{BKS'84}. However, the dependence with respect to the terminal time $T$ would not be explicit enough to infer an explicit observability constant.
\end{remark}

\section{The Carleman estimate}\label{sec:carleman}
In this section, we prove a boundary Carleman estimate for the homogeneous backward adjoint problem \eqref{eq:adj}, which is the key tool to prove the boundary null controllability of system \eqref{eq:intr}.

In the following lemma, we recall a weight function needed for the global boundary Carleman estimate. We refer to \cite[Lemma 2.1]{YY'09} for a precise statement and also for the proof.
\begin{lemma}[see \cite{Im'95}]\label{lem:eta} Given a nonempty relatively open set $\gamma \Subset\Gamma^0 \Subset  \Gamma$, there is a function 
$\eta\in C^2(\overline{\Omega})$ such that
\begin{equation}
\begin{aligned}
\mathrm{(i)} \hspace{0.6cm} \eta &> 0 &&\text{ in }\Omega,\\
\mathrm{(ii)} \; |\nabla \eta| &\geq c_0> 0 &&\text{ in }\overline{\Omega},\\
\mathrm{(iii)}\;\, \partial_\nu \eta &\leq -c_0 &&\text{ on }\Gamma \setminus \gamma,\\
\mathrm{(iv)}\quad\;\, \eta&=0, \; \nabla_{\Gamma}\eta =0 &&\text{ on }\Gamma\setminus \gamma.
\end{aligned}
\end{equation}
\end{lemma}

Let $\gamma\Subset \Gamma^0$ be a relatively open and non-empty, and let $\eta$ be as above. We then define the weight functions $\alpha$ and $\xi$ by
\begin{equation}
\begin{aligned}
\theta(t)&=\frac{1}{t(T-t)},  &&t\in (0,T),\\
\alpha(t,x)&=\theta(t) \left(\mathrm{e}^{2\lambda \|\eta\|_\infty}- \mathrm{e}^{\lambda \eta(x)}\right), &&(t,x) \in \overline{\Omega}_T,\\
\xi(t,x)&=\theta(t) \mathrm{e}^{\lambda \eta(x)},  &&(t,x) \in \overline{\Omega}_T,
\end{aligned}
\end{equation}
where $\lambda >1$ is a sufficiently large parameter. The weights $\alpha$ and $\xi$ are of class $C^2$, strictly positive on $\overline{\Omega}_T$, and blow up at $t=0$ and $t=T$. Moreover, we have the following formulas:
\begin{itemize}
\item On $\Gamma$, \quad $\nabla_{\Gamma} \alpha =-\lambda \xi \nabla_{\Gamma}\eta$ \;\;\; and \quad $\Delta_{\Gamma}\alpha=-\lambda \xi\Delta_{\Gamma}\eta-\lambda^2\xi|\nabla_{\Gamma}\eta|^2.$
 \item In $\overline{\Omega}_T$, \quad $\partial_t \alpha= (2t-T)\theta \alpha $\,\; and \quad $\partial_t^2 \alpha= 2\theta \alpha  +2(2t-T)\theta^2\alpha$.
\end{itemize}
Next, we state the Carleman estimate. In the proof, we follow the strategy of Lemma~3.2 of \cite{MMS'17}. In our setting, several new boundary terms arise from the new weight function, which is not null on the bundary as  the one of   \cite{MMS'17}, and these should be absorbed.
\begin{lemma}[Boundary Carleman estimate]\label{lem:carleman} 
There exist two positive constants $C=C\left(\Omega,\Gamma^0\right)$ and $\lambda_1=C\left(\Omega,\Gamma^0\right) \ge 1$ that are independent of $T$, and a positive constant $s_1=C\left(\Omega,\Gamma^0\right)\left(T + T^{8/3}\right)$ such that the strong solution $(\varphi,\varphi_\Gamma)\in \mathbb{E}_1$ of the backward system \eqref{eq:adj} satisfies
\begin{align} \label{Car}
& s^{-1} \int_{\Omega_T} \mathrm{e}^{-2s\alpha} \xi^{-1} \left(|\partial_t\varphi|^2 + |\Delta \varphi|^2\right) \, \d x\, \d t + s\lambda^2 \int_{\Omega_T} \mathrm{e}^{-2s\alpha} \xi |\nabla \varphi|^2 \, \d x\, \d t \notag\\
  & + s^{-1}\int_{\Gamma_T} \mathrm{e}^{-2s\alpha}\xi^{-1}\left(|\partial_t \varphi_\Gamma|^2 
   + |\Delta_\Gamma \varphi_\Gamma|^2\right) \, \d S\, \d t + s\lambda \int_{\Gamma_T} \mathrm{e}^{-2s\alpha} \xi |\nabla_\Gamma \varphi_\Gamma|^2 \, \d S\, \d t \notag\\
 &  + s^3\lambda^4\int_{\Omega_T} \mathrm{e}^{-2s\alpha} \xi^3 |\varphi|^2 \, \d x\, \d t +  s^3\lambda^3 \int_{\Gamma_T} \mathrm{e}^{-2s\alpha} \xi^3 |\varphi_\Gamma|^2 \, \d S\, \d t \notag\\
& + s\lambda \int_{\Gamma_T} \mathrm{e}^{-2s\alpha} \xi |\partial_\nu \varphi|^2 \, \d S\, \d t \notag\\
& \leq  C s^6\lambda^6 \int_{\Gamma^0_T} \mathrm{e}^{-s\alpha} \xi^{6} |\varphi_\Gamma|^2 \, \d S\, \d t
\end{align}
for all $\lambda \geq \lambda_1$, $s \geq s_1$.
\end{lemma}

In the following remark, we highlight a few facts about the boundary Carleman inequality \eqref{Car}.
\begin{remark}
It is worth to mention that:
\begin{itemize}
\item In the boundary Carleman inequality for the classical heat equation, the local term
\begin{equation}\label{lt}
s\lambda \int_{\Gamma^0_T} \mathrm{e}^{-2s\alpha} \xi (\partial_\nu \eta) |\partial_\nu \varphi|^2 \, \d S\, \d t 
\end{equation}
naturally appears in the right-hand side of the Carleman estimate, see, for instance, \cite[Théorème 3.4]{Ch'09}. A main difficulty in comparison with the classical boundary conditions is that, in our case, we must absorb the normal derivative term \eqref{lt} by the left-hand side of \eqref{Car}. This will be done by using the regularity estimates stated in Section \ref{sec:well}. The underlying term has been easily transmitted to the left-hand side of \eqref{Car} in the case of internal controllability, since in that case the weight function satisfies $\partial_\nu \eta^0 <0$ on whole $\Gamma$.
\item In contrast to the interior weight function used in \cite{MMS'17} for dynamic boundary conditions, the boundary weight function $\eta$ is not null on $\Gamma$, which makes the situation a bit more complicated by adding several new boundary terms in the proof of the Carleman estimate.
\end{itemize}
\end{remark}

\begin{proof}[Proof of Lemma \ref{lem:carleman}]
By a density argument, it suffices to consider smooth solutions $\varphi\in C^\infty\left([0,T]\times \overline{\Omega}\right)$. We will often write $\varphi$ instead of $\varphi_\Gamma$. Throughout the proof, $C$ will denote a generic constant which 
does not depend on $\lambda$, $s$ and $\varphi$.

We will also shorten some calculations done in \cite{MMS'17} and will pay a special attention to terms that need to be treated differently.
\bigskip

\noindent\textbf{Step 1. Change of variables.}
Let $\varphi\in C^\infty\left([0,T]\times \overline{\Omega}\right)$ satisfying \eqref{eq:adj},  $\lambda \geq \lambda_1\ge1$ and $s\ge s_1\ge 1$ be given. We consider the change of variable $\psi:= \mathrm{e}^{-s\alpha} \varphi$. Simple calculations yield that
\begin{align}
\nabla \alpha & = - \nabla \xi = - \lambda \xi \nabla \eta,\label{eq:xi2}\\
\Delta \alpha &
     = - \lambda^2 \xi |\nabla \eta|^2  - \lambda \xi \Delta \eta,\notag \\
\partial_t \psi & = \mathrm{e}^{-s\alpha} \partial_t \varphi - s \psi\partial_t\alpha, \notag\\
\nabla \psi & = \mathrm{e}^{-s\alpha} \nabla \varphi - s \psi\nabla\alpha 
 = \mathrm{e}^{-s\alpha} \nabla \varphi + s \lambda \psi \xi \nabla \eta, \label{eq:psi}\\
\Delta \psi &
   = \mathrm{e}^{-s\alpha} \Delta \varphi - 2s \nabla \psi \cdot \nabla \alpha - s^2 \psi |\nabla \alpha|^2 
     - s  \psi \Delta \alpha. \notag 
\end{align}
Then the equation fulfilled by $\psi$ in $\Omega_T$ is given by
\begin{align}
\partial_t \psi + d\Delta \psi & =- s\psi \partial_t\alpha + d s  \lambda^2 \psi \xi |\nabla \eta|^2  + ds \lambda \psi \xi \Delta \eta
  +  2 d s \lambda \xi \nabla \psi \cdot \nabla \eta \notag\\
  &\quad - d s^2\lambda^2  \psi \xi^2 |\nabla \eta |^2. \label{m1}
\end{align}
Similarly, on $\Gamma_T$ we obtain
\begin{align}
&\partial_t \psi + \delta \Delta_\Gamma \psi - d\partial_\nu \psi
  = - s\psi \partial_t\alpha-\delta s \psi\Delta_{\Gamma}\alpha -2\delta s \nabla_{\Gamma}\psi\cdot \nabla_{\Gamma} \alpha  +s\psi\partial_\nu \alpha \notag 
  \\&=- s\psi \partial_t\alpha +\delta s\lambda \xi\psi\Delta_{\Gamma}\eta+\delta s\lambda^2\xi\psi|\nabla_{\Gamma}\eta|^2+2\delta \lambda s\xi \nabla_{\Gamma}\psi\cdot \nabla_{\Gamma} \eta -\delta s^2\lambda^2 \xi^2\psi|\nabla_{\Gamma}\eta|^2\notag\\
  & \quad -ds\lambda \psi \xi  \partial_\nu \eta .\label{m2}
\end{align}
Following the decomposition in \cite{MMS'17}, we write the identities \eqref{m1} and \eqref{m2} as
\begin{equation}\label{m3}
 M_1\psi + M_2\psi = \tilde f \quad \text{in }\Omega_T, \qquad 
  N_1\psi + N_2\psi = \tilde g \quad \text{on }\Gamma_T,
\end{equation}
where
\begin{align*}
M_1\psi &= -2 d s  \lambda^2  \xi\psi |\nabla \eta|^2  - 2 d s \lambda \xi \nabla \psi  \cdot \nabla \eta 
           + \partial_t\psi = (M_1\psi)_1 + (M_1\psi)_2 + (M_1\psi)_3 , \\ 
M_2\psi & = d s^2\lambda^2 \xi^2\psi  |\nabla \eta |^2 + d \Delta \psi + s\psi \partial_t\alpha 
            = (M_2\psi)_1 + (M_2\psi)_2 + (M_2\psi)_3,\\ 
N_1\psi & = -2 \delta s\lambda^2\xi\psi|\nabla_{\Gamma}\eta|^2 -2 \delta \lambda s\xi \nabla_{\Gamma}\psi\cdot \nabla_{\Gamma} \eta + \partial_t\psi + d s\lambda \xi \psi  \partial_\nu \eta \\
     &= (N_1\psi)_1 + (N_1\psi)_2+ (N_1\psi)_3+ (N_1\psi)_4,\\ 
N_2\psi & = \delta s^2\lambda^2  \psi \xi^2 |\nabla_{\Gamma} \eta |^2 +\delta \Delta_\Gamma \psi + s\psi \partial_t\alpha - d \partial_\nu \psi 
     = (N_2\psi)_1 + (N_2\psi)_2 + (N_2\psi)_3 \notag\\
     & \quad +(N_2\psi)_4, \\
\tilde f & = d s \lambda \psi \xi \Delta \eta - d s \lambda^2 \psi \xi |\nabla \eta|^2 ,\\
\tilde g &= \delta s \lambda\psi\xi\Delta_{\Gamma}\eta - \delta s\lambda^2\xi\psi |\nabla_{\Gamma} \eta|^2 .
\end{align*}
Applying $\|\cdot\|_{L^2(\Omega_T)}^2$ resp.\ $\|\cdot\|_{L^2(\Gamma_T)}^2$ to the equations in \eqref{m3} 
and adding the resulting identities, we obtain
\begin{align}
& \|\tilde f\|_{L^2(\Omega_T)}^2 + \|\tilde g\|_{L^2(\Gamma_T)}^2
 = \|M_1\psi\|_{L^2(\Omega_T)}^2 + \|M_2\psi\|_{L^2(\Omega_T)}^2  +\|N_1\psi\|_{L^2(\Gamma_T)}^2 \notag\\
   & \quad + \|N_2\psi\|_{L^2(\Gamma_T)}^2 
 + 2 \sum_{i,j} \langle  (M_1\psi)_i, (M_2\psi)_j\rangle_{L^2(\Omega_T)} 
  + 2 \sum_{i,j} \langle (N_1\psi)_i, (N_2\psi)_j\rangle_{L^2( \Gamma_T)}.  \label{m12}
\end{align}
\smallskip

\noindent\textbf{Step 2. Estimating the mixed terms from below.}
First, we expand the mixed terms in \eqref{m12} by using, several times, integration by parts with respect to $t$, Green formula over $\Omega$ and the surface divergence theorem on $\Gamma$.

The first term is nonpositive
$$\langle (M_1\psi)_1, (M_2\psi)_1\rangle_{L^2(\Omega_T)} 
       = - 2 d^2 s^3\lambda^4 \int_{\Omega_T} |\nabla \eta|^4 \xi^3 \psi^2 \, \d x\, \d t.$$
Using integration by parts with \eqref{eq:xi2}, we derive
\begin{align*}
 & \hspace{2cm} \langle (M_1\psi)_2\, , (M_2\psi)_1\rangle_{ L^2(\Omega_T)} \\
& = 3 d^2 s^3 \lambda^4 \int_{\Omega_T} |\nabla \eta|^4 \xi^3 \psi^2 \, \d x\, \d t
    + d^2 s^3 \lambda^3 \int_{\Omega_T} \Delta \eta |\nabla \eta|^2 \xi^3 \psi^2 \, \d x\, \d t\\
& \quad + d^2 s^3\lambda^3\int_{\Omega_T}(\nabla (|\nabla \eta|^2)\cdot \nabla \eta)\xi^3\psi^2 \, \d x\, \d t
  -d^2 s^3\lambda^3 \int_{\Gamma_T} |\nabla \eta|^2 \xi^3  \partial_\nu \eta \psi^2  \, \d S\, \d t.
\end{align*}
Integrating by parts in $t$ and using $\psi(0)=\psi(T)=0$, we obtain
\begin{align*}
 \langle (M_1\psi)_3, (M_2\psi)_1\rangle_{L^2(\Omega_T)}
   &= -d s^2\lambda^2 \int_{\Omega_T} |\nabla\eta|^2  \partial_t \xi \xi \psi^2 \, \d x\, \d t.
\end{align*}
Integration by parts yields
\begin{align*}
& \hspace{2cm} \langle (M_1\psi)_1, (M_2\psi)_2\rangle_{L^2(\Omega_T)} \\
& = 2 d^2 s\lambda^2 \int_{\Omega_T} |\nabla \eta|^2 \xi   |\nabla \psi|^2  \, \d x\, \d t
    + 2 d^2 s\lambda^2 \int_{\Omega_T} \xi \psi \nabla (|\nabla \eta|^2)\cdot  \nabla \psi  \, \d x\, \d t \\
&\qquad + 2 d^2 s\lambda^3 \int_{\Omega_T}|\nabla \eta|^2 \psi  \xi \nabla \eta \cdot\nabla\psi  \, \d x\, \d t
 -  2 d^2 s\lambda^2 \int_{\Gamma_T} |\nabla \eta|^2 \xi \psi \partial_\nu \psi \, \d S\, \d t.
\end{align*}
The next summand is given by 
\begin{align}
& \hspace{2cm} \langle (M_1\psi)_2, (M_2\psi)_2\rangle_{L^2(\Omega_T)} \notag\\
& = -2 d^2 s\lambda \int_{\Gamma_T}  \nabla \eta \cdot \nabla \psi   \partial_\nu \psi\xi \, \d S\, \d t 
    + 2 d^2 s\lambda \int_{\Omega_T} \xi\, (\nabla^2 \eta \,\nabla\psi) \cdot \nabla \psi \, \d x\, \d t \label{D3}\\
& \qquad    + 2 d^2 s\lambda^2 \int_{\Omega_T}  \xi |\nabla \eta \cdot \nabla\psi|^2 \, \d x\, \d t
   + d^2 s\lambda \int_{\Omega_T}  \nabla \eta \cdot \nabla|\nabla\psi|^2\,\xi  \, \d x\, \d t, \notag
\end{align}
using integration by parts and \eqref{eq:xi2}. Here, $\nabla^2 \eta$ stands for the standard Hessian matrix of the function $\eta$. For the last term in \eqref{D3}, integration by parts and \eqref{eq:xi2} yield
\begin{align*}
& \hspace{1cm} d^2 s\lambda \int_{\Omega_T}  \nabla \eta \cdot \nabla|\nabla\psi|^2\,\xi  \, \d x\, \d t \, \\
& = d^2 s\lambda \int_{\Gamma_T}  \partial_\nu  \eta |\nabla \psi|^2\,\xi  \, \d S\, \d t
  - d^2 s\lambda^2 \int_{\Omega_T}   \xi   |\nabla \eta|^2 |\nabla\psi|^2  \, \d x\, \d t \\
& \quad \; - d^2 s\lambda \int_{\Omega_T}   \xi \Delta  \eta |\nabla\psi|^2  \, \d x\, \d t.
\end{align*}
Next,
\begin{align*}
& \hspace{1cm} -2 d^2 s\lambda \int_{\Gamma_T}  \xi\nabla \eta \cdot \nabla \psi   \partial_\nu \psi \, \d S\, \d t \\
&= -2 d^2 s\lambda \int_{\Gamma_T}  \xi(\nabla_\Gamma \eta + (\partial_\nu \eta)\nu) \cdot (\nabla_\Gamma \psi + (\partial_\nu \psi)\nu) \partial_\nu \psi \, \d S\, \d t \\
& =-2 d^2 s\lambda \int_{\Gamma_T}  \xi\partial_\nu \eta (\partial_\nu \psi)^2 \, \d S\, \d t -2 d^2 s\lambda \int_{\Gamma_T} \xi \partial_\nu \psi(\nabla_\Gamma \eta \cdot \nabla_\Gamma \psi ) \, \d S\, \d t.
\end{align*}
Since $\nabla \psi$ vanishes at $t = 0$ and $t = T$ in view of \eqref{eq:psi}, we obtain
\begin{align}
 \langle (M_1\psi)_3, (M_2\psi)_2\rangle_{L^2(\Omega_T)}
 & = d \int_{\Gamma_T} \partial_t \psi \partial_\nu \psi \, \d S\, \d t. \label{m13} \\
\langle (M_1\psi)_1, (M_2\psi)_3\rangle_{L^2(\Omega_T)} 
 &= - 2 d s^2\lambda^2 \int_{\Omega_T} |\nabla \eta|^2 \partial_t \alpha \xi \psi^2 \, \d x\, \d t. \notag
\end{align}
Integration by parts and \eqref{eq:xi2} imply
\begin{align*}
& \hspace{2cm} \langle (M_1\psi)_2, (M_2\psi)_3\rangle_{L^2(\Omega_T)} \\
& = - d s^2 \lambda \int_{\Gamma_T} \partial_t \alpha \xi \partial_\nu\eta \psi^2 \, \d S\, \d t
  + d s^2 \lambda \int_{\Omega_T} \nabla (\partial_t \alpha) \cdot \nabla\eta \xi \psi^2 \, \d x\, \d t\\
&\qquad + d s^2 \lambda^2 \int_{\Omega_T} \partial_t \alpha \xi  |\nabla\eta|^2 \psi^2 \, \d x\, \d t
   + d s^2 \lambda \int_{\Omega_T} \partial_t \alpha \xi \Delta \eta \psi^2 \, \d x\, \d t.
\end{align*}
Integrating by parts with respect to $t$, we can derive
\begin{equation}\label{eq:M33}
\langle (M_1\psi)_3, (M_2\psi)_3\rangle_{L^2(\Omega_T)} 
  =  - \frac{s}{2} \int_{\Omega_T} \partial_t^2 \alpha \psi^2\, \d x\, \d t,
\end{equation}
since $\psi$ vanishes at the endpoints.

We now consider the boundary terms $N_1$ and $N_2$, employing the surface divergence  theorem \eqref{eq:sdt} several times, 
we have
$$\langle (N_1\psi)_1, (N_2\psi)_1\rangle_{L^2(\Gamma_T)}=- 2 \delta^2 s^3 \lambda^4 \int_{\Gamma_T} |\nabla_{\Gamma} \eta|^4 \xi^3 \psi^2 \, \d S\, \d t.$$
Using integration by parts, we further derive
\begin{align*}
\langle (N_1\psi)_2&\, , (N_2\psi)_1\rangle_{ L^2(\Gamma_T)} 
= - \delta^2 s^3\lambda^3 \int_{\Gamma_T} |\nabla_{\Gamma} \eta|^2 \xi^3  \nabla_{\Gamma} \eta \cdot \nabla_{\Gamma} (\psi^2) \, \d S\, \d t\\
&= \delta^2 s^3\lambda^3 \int_{\Gamma_T} \text{div}_{\Gamma}\big(|\nabla_\Gamma \eta|^2 \xi^3\nabla_\Gamma \eta\big)\psi^2 \, \d S\, \d t\\
& = 3 \delta^2 s^3 \lambda^4 \int_{\Gamma_T} |\nabla_{\Gamma} \eta|^4 \xi^3 \psi^2 \, \d S\, \d t
+ \delta^2 s^3 \lambda^3 \int_{\Gamma_T} \Delta_{\Gamma} \eta |\nabla_{\Gamma} \eta|^2 \xi^3 \psi^2 \, \d S\, \d t\\
& \qquad + \delta^2 s^3\lambda^3\int_{\Gamma_T}(\nabla_{\Gamma} (|\nabla_{\Gamma} \eta|^2)\cdot \nabla_{\Gamma} \eta)\xi^3\psi^2 \, \d S\, \d t.
\end{align*}
Therefore
\begin{align}\label{eq:N1}
\langle (N_1\psi)_1&, (N_2\psi)_1\rangle_{L^2(\Gamma_T)}  + \langle (N_1\psi)_2 , 
(N_2\psi)_1\rangle_{L^2(\Gamma_T)}\notag\\
& = \delta^2 s^3 \lambda^4 \int_{\Gamma_T} |\nabla_{\Gamma} \eta|^4 \xi^3 \psi^2 \, \d S\, \d t
+ \delta^2 s^3 \lambda^3 \int_{\Gamma_T} \Delta_{\Gamma} \eta |\nabla_{\Gamma} \eta|^2 \xi^3 \psi^2 \, \d S\, \d t\notag\\
& \qquad + \delta^2 s^3\lambda^3\int_{\Gamma_T}(\nabla_{\Gamma} (|\nabla_{\Gamma} \eta|^2)\cdot \nabla_{\Gamma} \eta)\xi^3\psi^2 \, \d S\, \d t.
\end{align}

\begin{align*}
\langle (N_1\psi)_3, (N_2\psi)_1\rangle_{L^2(\Gamma_T)} 
&= \frac{\delta}{2} s^2\lambda^2 \int_{\Gamma_T} \partial_{t}(\psi^2) \xi^2 |\nabla_{\Gamma} \eta|^2 \, \d S\, \d t\\
&= - \delta s^2\lambda^2 \int_{\Gamma_T} |\nabla_{\Gamma} \eta|^2 \partial_{t} \xi \xi \psi^2 \, \d S\, \d t.
\end{align*}
The next terms are given by  
\begin{align*}
\langle (N_1\psi)_4, (N_2\psi)_1\rangle_{L^2(\Gamma_T)} 
&= d\delta s^3\lambda^3 \int_{\Gamma_T}  \xi^3  \partial_\nu \eta |\nabla_{\Gamma} \eta|^2 \psi^2 \, \d S\, \d t.
\end{align*}

\begin{align*}
\langle (N_1\psi)_1, (N_2\psi)_2\rangle_{L^2(\Gamma_T)} &= -2 \delta^2 s\lambda^2 \int_{\Gamma_T} |\nabla_\Gamma \eta|^2 \xi \psi \Delta_\Gamma \psi \, \d S\, \d t\\
& = 2\delta^2 s\lambda^2\int_{\Gamma_T}\nabla_{\Gamma}(\xi \psi |\nabla_{\Gamma}\eta|^2)\cdot \nabla_{\Gamma}\psi \, \d S\, \d t\\
&\hspace{-3.5cm} =2\delta^2 s\lambda^2 \int_{\Gamma_T} |\nabla_{\Gamma} \eta|^2 \xi |\nabla_{\Gamma} \psi|^2 \, \d S\, \d t + 2\delta^2 s\lambda^2\int_{\Gamma_T}\xi\psi \nabla_{\Gamma}(|\nabla_{\Gamma} \eta|^2)\cdot \nabla_{\Gamma}\psi \, \d S\, \d t \\
& \hspace{-3cm} + 2\delta^2 s\lambda^3 \int_{\Gamma_T} |\nabla_{\Gamma}\eta|^2 \psi\xi (\nabla_{\Gamma}\eta\cdot\nabla_{\Gamma}\psi ) \, \d S\, \d t.
\end{align*}

\begin{align}
\langle (N_1\psi)_2&, (N_2\psi)_2\rangle_{L^2(\Gamma_T)} 
= -2\delta^2 s\lambda \int_{\Gamma_T} (\nabla_{\Gamma} \eta \cdot \nabla_{\Gamma}\psi) \xi  \Delta_{\Gamma} \psi \, \d S\, \d t \notag\\
& = 2\delta^2 s\lambda\int_{\Gamma_T}\nabla_{\Gamma}\big((\nabla_{\Gamma} \eta \cdot \nabla_{\Gamma}\psi)\xi \big) \cdot \nabla_{\Gamma} \psi \, \d S\, \d t\notag\\
& =2\delta^2 s\lambda \int_{\Gamma_T} \xi\, (\nabla_{\Gamma}^2 \eta \,\nabla_{\Gamma}\psi) \cdot \nabla_{\Gamma} \psi \, \d x\, \d t 
+ 2\delta^2 s\lambda^2 \int_{\Gamma_T}  \xi |\nabla_{\Gamma} \eta \cdot \nabla_{\Gamma}\psi|^2 \, \d S\, \d t \notag\\
& \qquad + \delta^2 s\lambda \int_{\Gamma_T}  \xi \nabla_{\Gamma} \eta \cdot \nabla_{\Gamma}|\nabla_{\Gamma}\psi|^2 \, \d S\, \d t \notag\\
& =2\delta^2 s\lambda \int_{\Gamma_T} \xi\, (\nabla_{\Gamma}^2 \eta \,\nabla_{\Gamma}\psi) \cdot \nabla_{\Gamma} \psi \, \d x\, \d t
+ 2\delta^2 s\lambda^2 \int_{\Gamma_T}  \xi |\nabla_{\Gamma} \eta \cdot \nabla_{\Gamma}\psi|^2 \, \d S\, \d t \notag\\
& \qquad - \delta^2 s\lambda \int_{\Gamma_T} \mathrm{div}_\Gamma(\xi \nabla_{\Gamma} \eta) |\nabla_{\Gamma}\psi|^2 \, \d S\, \d t \notag\\
&= 2\delta^2 s\lambda \int_{\Gamma_T} \xi\, (\nabla_{\Gamma}^2 \eta \,\nabla_{\Gamma}\psi) \cdot \nabla_{\Gamma} \psi \, \d x\, \d t
+ 2\delta^2 s\lambda^2 \int_{\Gamma_T}  \xi |\nabla_{\Gamma} \eta \cdot \nabla_{\Gamma}\psi|^2 \, \d S\, \d t \notag\\
& \qquad  - \delta^2 s\lambda^2 \int_{\Gamma_T}  \xi |\nabla_{\Gamma} \eta|^2 |\nabla_{\Gamma}\psi|^2 \, \d S\, \d t - \delta^2 s\lambda \int_{\Gamma_T}  \xi \Delta_\Gamma \eta |\nabla_{\Gamma}\psi|^2 \, \d S\, \d t.
\end{align}
Since $\psi(0)=\psi(T)=0$, we have
$$\langle (N_1\psi)_3, (N_2\psi)_2\rangle_{L^2(\Gamma_T)} 
 = 0.$$
Moreover, using   $\nabla_\Gamma \xi =\lambda \xi \nabla_\Gamma \eta$, one has
\begin{align*}
& \hspace{2cm} \langle (N_1\psi)_4 , (N_2\psi)_2\rangle_{L^2(\Gamma_T)} 
   \\
& = - d \delta s\lambda \int_{\Gamma_T}  \xi \psi \langle \nabla_\Gamma\partial_\nu \eta , 
    \nabla_\Gamma \psi\rangle_\Gamma \, \d S\, \d t 
    - d \delta s\lambda \int_{\Gamma_T} \partial_\nu \eta \xi  |\nabla_\Gamma \psi|^2 \, \d S\, \d t\\
    &\quad\;  - d \delta s\lambda^2 \int_{\Gamma_T} \partial_\nu \eta \xi \psi \nabla_\Gamma \eta \cdot \nabla_\Gamma \psi \, \d S\, \d t.
\end{align*}
$$\langle (N_1\psi)_1, (N_2\psi)_3\rangle_{L^2(\Gamma_T)} = -2 \delta s^2\lambda^2  \int_{\Gamma_T} |\nabla_{\Gamma} \eta |^2\partial_t \alpha \xi \psi^2 \, \d S\, \d t.$$
Integration by parts and \eqref{eq:xi2} imply
\begin{align*}
&\langle (N_1\psi)_2, (N_2\psi)_3\rangle_{L^2(\Gamma_T)} 
= - \delta s^2 \lambda \int_{\Gamma_T} \partial_t \alpha \xi \nabla_\Gamma\eta \cdot \nabla_\Gamma (\psi^2) \, \d S\, \d t\\
&= \delta s^2 \lambda \int_{\Gamma_T} \text{div}_{\Gamma}(\partial_t \alpha \xi \nabla_{\Gamma}\eta) \psi^2 \, \d S\, \d t\\
& = \delta s^2 \lambda \int_{\Gamma_T} \nabla_{\Gamma} (\partial_t \alpha) \cdot \nabla_{\Gamma}\eta \xi \psi^2 \, \d S\, \d t\\
&\qquad + \delta s^2 \lambda^2 \int_{\Gamma_T} \partial_t \alpha \xi |\nabla_{\Gamma}\eta|^2 \psi^2 \, \d S\, \d t
+ \delta s^2 \lambda \int_{\Gamma_T} \partial_t \alpha \xi \Delta_{\Gamma} \eta \psi^2 \, \d S\, \d t.
\end{align*}
Hence
\begin{align*}
\langle (N_1\psi)_3, (N_2\psi)_3\rangle_{L^2(\Gamma_T)} &= -\frac{s}{2} \int_{\Gamma_T} \partial_t^2 \alpha (\psi^2) \, \d S\, \d t. \\
 \langle (N_1\psi)_4, (N_2\psi)_3\rangle_{L^2(\Gamma_T)} &= ds^2\lambda \int_{\Gamma_T} \partial_\nu \eta \partial_t \alpha \xi \psi^2 \, \d S\, \d t.\\
\langle (N_1\psi\,)_1, (N_2\psi)_4\rangle_{L^2(\Gamma_T)} 
     &=2 d \delta s\lambda^2\int_{\Gamma_T}\xi\psi\partial_{\nu}\psi|\nabla_{\Gamma}\eta|^2 \, \d S\, \d t.\\
     \langle (N_1\psi\,)_2, (N_2\psi)_4\rangle_{L^2(\Gamma_T)}
     &=2 d \delta s\lambda \int_{\Gamma_T}\xi\partial_{\nu}\psi (\nabla_{\Gamma}\psi \cdot \nabla_{\Gamma} \eta ) \, \d S\, \d t. \\
 \langle (N_1\psi\,)_3, (N_2\psi)_4\rangle_{L^2(\Gamma_T)} 
     &= -d \int_{\Gamma_T} \partial_t \psi \partial_\nu \psi \, \d S\, \d t.
\end{align*}
The last term cancels with the one from \eqref{m13}. Next,
\begin{align}
\langle (N_1\psi)_4, (N_2\psi)_4\rangle_{L^2(\Gamma_T)} = - d^2 s\lambda \int_{\Gamma_T}  \xi \partial_\nu \eta \partial_\nu \psi \psi \, \d S\, \d t.
\end{align}
Summing up all terms together, we derive
\begin{align*}
& \sum_{i,j} \langle (M_1\psi)_i, (M_2\psi)_j\rangle_{L^2(\Omega_T)}
\\
& = d^2 s^3 \lambda^4 \int_{\Omega_T} |\nabla \eta|^4 \xi^3 \psi^2 \, \d x\, \d t
+\ d^2 s^3 \lambda^3 \int_{\Omega_T} \Delta \eta |\nabla \eta|^2 \xi^3 \psi^2 \, \d x\, \d t\\
& \quad +\  d^2  s^3\lambda^3\int_{\Omega_T}(\nabla (|\nabla \eta|^2)\cdot \nabla \eta)\xi^3\psi^2 \, \d x\, \d t
-  d  s^2\lambda^2 \int_{\Omega_T} |\nabla\eta|^2  \partial_t \xi\xi\psi^2 \, \d x\, \d t\\
& \quad -  d  s^2\lambda^2 \int_{\Omega_T} |\nabla \eta|^2 \partial_t \alpha\xi\psi^2 \, \d x\, \d t +\  d s^2 \lambda \int_{\Omega_T} \nabla (\partial_t \alpha) \cdot \nabla\eta \xi\psi^2 \, \d x\, \d t \\
& \quad +\  d  s^2 \lambda \int_{\Omega_T} \partial_t \alpha\xi\Delta \eta \psi^2 \, \d x\, \d t - \frac{s}{2} \int_{\Omega_T} \partial_t^2 \alpha\psi^2 \, \d x\, \d t\\
& \quad  +\ d^2 s\lambda^2 \int_{\Omega_T} \xi |\nabla \eta|^2  |\nabla \psi|^2  \, \d x\, \d t 
+ 2  d^2  s\lambda^2 \int_{\Omega_T} \xi\psi \nabla (|\nabla \eta|^2)\cdot  \nabla \psi  \, \d x\, \d t \\
& \quad + 2 d^2 s \lambda^3 \int_{\Omega_{T}} |\nabla \eta|^2 \psi \xi \nabla \eta \cdot \nabla \psi \, \d x\, \d t
+\ 2  d^2  s\lambda \int_{\Omega_T} \xi\, (\nabla^2 \eta \,\nabla\psi) \cdot \nabla \psi \, \d x\, \d t \\
& \quad -\  d^2  s\lambda \int_{\Omega_T}   \xi \Delta  \eta |\nabla\psi|^2  \, \d x\, \d t +\ 2  d^2  s\lambda^2 \int_{\Omega_T}  \xi|\nabla \eta \cdot \nabla\psi|^2 \, \d x\, \d t\\
& \quad - d^2 s^3\lambda^3 \int_{\Gamma_T} |\nabla \eta|^2 \xi^3  \partial_\nu \eta \psi^2  \, \d S\, \d t  -  2 d^2 s\lambda^2 \int_{\Gamma_T} |\nabla \eta|^2 \xi \psi \partial_\nu \psi \, \d S\, \d t \\
& \quad - d^2 s\lambda \int_{\Gamma_T}  \xi \partial_\nu \eta (\partial_\nu \psi)^2 \, \d S\, \d t   +  d^2 s\lambda \int_{\Gamma_T}  \partial_\nu  \eta |\nabla_\Gamma \psi|^2\,\xi  \, \d S\, \d t\\
& \quad -2 d^2 s\lambda \int_{\Gamma_T} \xi \partial_\nu \psi(\nabla_\Gamma \eta \cdot \nabla_\Gamma \psi )\, \d S\, \d t 
+ d \int_{\Gamma_T} \partial_t \psi \partial_\nu \psi \, \d S \, \d t \\
& \quad - d s^2 \lambda \int_{\Gamma_T} \partial_t \alpha \xi \partial_\nu\eta \psi^2 \, \d S\, \d t \\
& = \sum_{k=1}^{14} I_k + \sum_{k=1}^7 B_{1,k}.\\
& \sum_{i,j} \langle (N_1\psi)_i, (N_2\psi)_j\rangle_{L^2(\Gamma_T)} \\
& = \delta^2 s^3 \lambda^4 \int_{\Gamma_T} |\nabla_{\Gamma} \eta|^4 \xi^3 \psi^2 \, \d S\, \d t  + \delta^2 s^3 \lambda^3 \int_{\Gamma_T} \Delta_{\Gamma} \eta |\nabla_{\Gamma} \eta|^2 \xi^3 \psi^2 \, \d S\, \d t\\
& \quad + \delta^2 s^3\lambda^3\int_{\Gamma_T}(\nabla_{\Gamma} (|\nabla_{\Gamma} \eta|^2)\cdot \nabla_{\Gamma} \eta)\xi^3\psi^2 \, \d S\, \d t  - \delta s^2\lambda^2 \int_{\Gamma_T} |\nabla_{\Gamma} \eta|^2 \partial_{t} \xi \xi \psi^2 \, \d S\, \d t \\
& \quad - \delta s^2\lambda^2  \int_{\Gamma_T} |\nabla_{\Gamma} \eta |^2\partial_t \alpha \xi \psi^2 \, \d S\, \d t  + \delta s^2 \lambda \int_{\Gamma_T} \nabla_{\Gamma} (\partial_t \alpha) \cdot \nabla_{\Gamma}\eta \xi \psi^2 \, \d S\, \d t \\
& \quad + \delta s^2 \lambda \int_{\Gamma_T} \partial_t \alpha \xi \Delta_{\Gamma} \eta \psi^2 \, \d S\, \d t -\frac{s}{2} \int_{\Gamma_T} \partial_t^2 \alpha \psi^2 \, \d S\, \d t\\
& \quad  +\ \delta^2 s\lambda^2 \int_{\Gamma_T} \xi |\nabla_{\Gamma} \eta|^2 |\nabla_{\Gamma} \psi|^2 \, \d S\, \d t  + 2\delta^2 s\lambda^2\int_{\Gamma_T}\xi\psi \nabla_{\Gamma}(|\nabla_{\Gamma} \eta|^2)\cdot \nabla_{\Gamma}\psi \, \d S\, \d t \\
&  +\ 2 \delta^2 s \lambda^3 \int_{\Gamma_T} |\nabla_{\Gamma} \eta|^2 \psi \xi (\nabla_{\Gamma} \eta \cdot \nabla_{\Gamma} \psi) \, \d S\, \d t
 + 2\delta^2 s\lambda \int_{\Gamma_T} \xi\, (\nabla_{\Gamma}^2 \eta \,\nabla_{\Gamma}\psi) \cdot \nabla_{\Gamma} \psi \, \d S\, \d t \\
& \quad  - \delta^2 s\lambda \int_{\Gamma_T}  \xi \Delta_\Gamma \eta |\nabla_{\Gamma}\psi|^2 \, \d S\, \d t +\ 2\delta^2 s\lambda^2 \int_{\Gamma_T}  \xi |\nabla_{\Gamma} \eta \cdot \nabla_{\Gamma}\psi|^2 \, \d S\, \d t \\
& \quad +  d\delta s^3\lambda^3 \int_{\Gamma_T}  \xi^3  \partial_\nu \eta |\nabla_{\Gamma} \eta|^2 \psi^2 \, \d S\, \d t   +   2 d \delta s\lambda^2\int_{\Gamma_T}\xi\psi\partial_{\nu}\psi|\nabla_{\Gamma}\eta|^2 \, \d S\, \d t \\
& \quad -  d \delta s\lambda \int_{\Gamma_T} \partial_\nu \eta \xi  |\nabla_\Gamma \psi|^2 \, \d S\, \d t   +   2 d \delta s\lambda \int_{\Gamma_T}\xi\partial_{\nu}\psi (\nabla_{\Gamma}\psi \cdot \nabla_{\Gamma} \eta ) \, \d S\, \d t \\
& \quad - d \int_{\Gamma_T} \partial_t \psi \partial_\nu \psi \, \d S \, \d t  + d s^2 \lambda \int_{\Gamma_T} \partial_t \alpha \xi \partial_\nu\eta \psi^2 \, \d S\, \d t \\
& \quad - d \delta s\lambda \int_{\Gamma_T} \xi \psi \langle \nabla_\Gamma(\partial_\nu \eta), \nabla_\Gamma \psi\rangle_\Gamma \, \d S\, \d t - d \delta s\lambda^2 \int_{\Gamma_T} \partial_\nu \eta \xi \psi \nabla_\Gamma \eta \cdot \nabla_\Gamma \psi \, \d S\, \d t \\
& \quad + 2\delta^2 s\lambda^3 \int_{\Gamma_T} |\nabla_{\Gamma}\eta|^2 \psi\xi (\nabla_{\Gamma}\eta\cdot\nabla_{\Gamma}\psi ) \, \d S\, \d t  -    d^2 s\lambda \int_{\Gamma_T}  \xi \partial_\nu \eta \partial_\nu \psi \psi\, \d S\, \d t   \\
& = \sum_{k=1}^{14} J_k + \sum_{k=1}^6 B_{2,k} + \sum_{k=1}^4 \widetilde{B}_{k}.
\end{align*}
Next, we estimate the resulting terms from below. We will often use the following basic estimates on $\overline{\Omega}_T$, 
\begin{equation}
\xi \le C T^2 \xi^2, \quad |\nabla \alpha| \leq C \lambda \xi, \quad |\partial_t \alpha| + |\partial_t \xi|\leq C T\xi^2, \quad |\partial_t^2\alpha| \le C T^2\xi^3. \label{m6}
\end{equation}
Using $|\nabla \eta| \ge c_0 >0$ in $\overline{\Omega}$ and $|\partial_\nu \eta| \ge c_0>0$ on $\Gamma \setminus \gamma$, we obtain
\begin{align}
\sum_{k=1}^{3} I_k + B_{1,1} &\ge C s^3\lambda^4 \int_{\Omega_T} \xi^3 \psi^2 \, \d x\, \d t + C s^3\lambda^3 \int_{(\Gamma \setminus \gamma)\times (0,T)} \xi^3 \psi^2 \, \d S\, \d t \notag\\
& \quad -C s^3\lambda^3 \int_{\gamma_T} \xi^3 \psi^2 \, \d S\, \d t
\end{align}
for $\lambda_1$ large enough. Since $J_1$ is nonnegative and $\nabla_\Gamma \eta =0$ on $\Gamma \setminus \gamma$, we infer that
\begin{align}
\sum_{k=1}^{3} J_k  + B_{2,1} \ge - C s^3\lambda^3 \int_{\gamma_T} \xi^3 \psi^2 \, \d S\, \d t.
\end{align}
Since $|\partial_t \xi | \le C T\xi^2$ and $|\partial_t^2 \alpha | \le C T^2 \xi^3$, we obtain
\begin{align}
\sum_{k=4}^{8} I_k  + \sum_{k=4}^{8} J_k  \ge - C s^2\lambda^2 \int_{\Omega_T} \xi^3 \psi^2 \, \d x\, \d t - C s^2\lambda^2 \int_{\Gamma_T} \xi^3 \psi^2 \, \d S\, \d t
\end{align}
Using again the fact $|\nabla \eta| \ge c_0 >0$ in $\overline{\Omega}$ and $|\partial_\nu \eta| \ge c_0>0$ on $\Gamma \setminus \gamma$, we derive
\begin{align}
I_9  + B_{1,3} + B_{2,3}  &\ge C s\lambda^2 \int_{\Omega_T} \xi |\nabla \psi|^2  \, \d x\, \d t  + C s\lambda \int_{(\Gamma \setminus \gamma)\times (0,T)} \xi  (\partial_\nu \psi)^2 \, \d S\, \d t \notag\\
& \quad - C s\lambda \int_{\gamma_T} \xi  (\partial_\nu \psi)^2 \, \d S\, \d t + C s\lambda \int_{(\Gamma \setminus \gamma)\times (0,T)} \xi  |\nabla_\Gamma \psi|^2 \, \d S\, \d t\notag\\
&  \quad - C s\lambda \int_{\gamma_T} \xi  |\nabla_\Gamma \psi|^2 \, \d S\, \d t.
\end{align}
Now, applying Young's inequality to 
 $(\lambda^2 \xi^{1/2}\psi) \big( \xi^{1/2}\nabla (|\nabla \eta|^2) \cdot \nabla \psi\big)$ together with
 $\big (s\lambda \psi  \xi  \big) \big( \nabla \eta  \cdot \nabla \psi\big)$, respectively, we find
\begin{align*}
\Big|2s\lambda^2 \int_{\Omega_T} \xi \psi \nabla (|\nabla \eta|^2)\cdot  \nabla \psi  \, \d x\, \d t \Big| 
&\leq C s \lambda^4 \int_{\Omega_T} \xi \psi^2  \, \d x\, \d t
  + C s \int_{\Omega_T} \xi |\nabla \psi|^2  \, \d x\, \d t,\\
\Big|2s\lambda^3 \int_{\Omega_T}|\nabla \eta|^2 \psi  \xi \nabla \eta \cdot\nabla\psi  \, \d x\, \d t \Big| 
&\leq C s^2\lambda^4 \int_{\Omega_T} \xi^2 \psi^2 \, \d x\, \d t 
     + C \lambda^2 \int_{\Omega_T} |\nabla \psi|^2 \, \d x\, \d t.
\end{align*}     
It follows
\begin{align}
I_{10} + I_{11} \ge - C s^2\lambda^4 \int_{\Omega_T} \xi^2\psi^2 \, \d x\, \d t - C\int_{\Omega_T}(s\xi +\lambda^2) |\nabla\psi|^2 \, \d x\, \d t.
\end{align}
By virtue of $\nabla_\Gamma \eta =0$ on $\Gamma \setminus \gamma$, it turns out that
\begin{align}
J_{10} + J_{11} &\ge - C s^2\lambda^4 \int_{\gamma_T} \xi^2 \psi^2 \, \d S\, \d t - C s\int_{\Gamma_T} \xi |\nabla_\Gamma \psi|^2 \, \d S\, \d t \notag\\
&\quad - C \lambda^2 \int_{\gamma_T} \xi |\nabla_\Gamma \psi|^2 \, \d S\, \d t.
\end{align}
Since $\eta \in C^2(\overline{\Omega})$ and $\nabla_\Gamma \eta =0$ on $\Gamma \setminus \gamma$, we can estimate
\begin{align*}
I_{12} + I_{13} + J_{12} + J_{13} \ge - C s\lambda \int_{\Omega_T} \xi |\nabla \psi|^2 \, \d x\, \d t - C s\lambda \int_{\gamma_T} \xi |\nabla_\Gamma \psi|^2 \, \d S\, \d t.
\end{align*}
The next terms are nonnegative
\begin{align}
I_{14} + J_9 + J_{14} \ge 0.
\end{align}
Young's inequality implies
\begin{align}
B_{1,2} + B_{2,2} + \widetilde{B}_{4} &\ge -Cs\int_{\Gamma_T}|\lambda^{3/2}\xi^{1/2}\psi|\,|\lambda^{1/2}\xi^{1/2}\partial_\nu \psi| \, \d S\, \d t \notag\\
& \ge -C\left( s^2 \lambda^3 \int_{\Gamma_T} \xi^3 \psi^2 \, \d S\, \d t 
    + \lambda \int_{\Gamma_T} \xi (\partial_\nu \psi)^2 \, \d S\, \d t\right).
\end{align}
Using the fact that $\xi(t,\cdot)$ is constant on $\Gamma \setminus \gamma$ (we denote it by $\xi(t)$) with the divergence formula, we then obtain
\begin{align}
& |B_{1,4}| \leq Cs\lambda \int_0^T \xi(t) \int_{\Gamma\setminus \gamma}  |\nabla_\Gamma \psi|^2 \, \d S\, \d t + Cs \lambda \int_{\gamma_T} \xi |\nabla_\Gamma \psi|^2 \, \d S\, \d t \notag \\
&\leq C \int_0^T \int_{\Gamma}\big(s^{-1/2} \xi^{-1/2}(t)  |\Delta_\Gamma \psi|\big) \big(s^{3/2} \lambda 
                   \xi^{3/2}(t) |\psi|\big) \, \d S\, \d t  + Cs \lambda \int_{\gamma_T} \xi |\nabla_\Gamma \psi|^2 \, \d S\, \d t\notag \\
& \leq  s^{-1} \int_{\Gamma_T} \xi^{-1} |\Delta_\Gamma \psi|^2 \, \d S\, \d t 
  + C s^3 \lambda^2 \int_{\Gamma_T} \xi^3 \psi^2 \, \d S\, \d t  +Cs \lambda \int_{\gamma_T} \xi |\nabla_\Gamma \psi|^2 \, \d S\, \d t.\notag 
\end{align}
We apply Young's inequality to $(s^{1/2}\xi^{1/2} |\partial_\nu \psi|)$ and $(s^{1/2} \lambda \xi^{1/2} |\nabla_\Gamma \eta| |\nabla_\Gamma \psi|)$ and we use $\nabla_\Gamma \eta =0$ on $\Gamma \setminus \gamma$ to obtain
\begin{align}
|B_{1,5} + B_{2,4}| & \le C s \lambda \int_{\Gamma_T} \xi |\partial_\nu \psi| |\nabla_\Gamma \eta| |\nabla_\Gamma \psi| \, \d S\, \d t \notag\\
& \le C s \int_{\Gamma_T} \xi (\partial_\nu \psi)^2 \, \d S\, \d t + C s \lambda^2 \int_{\gamma_T} \xi |\nabla_\Gamma \psi|^2 \, \d S\, \d t.
\end{align}
Moreover, we have
$$B_{1,6}+B_{2,5}+B_{1,7}+B_{2,6}=0.$$
Finally, 
\begin{align}
|\widetilde{B}_{1}| & \le C \left(s^2 \lambda^2 \int_{\Gamma_T} \xi \psi^2 \, \d S\, \d t + \int_{\Gamma_T} \xi  |\nabla_\Gamma \psi|^2 \, \d S\, \d t \right),\\
|\widetilde{B}_{2}| & \le C \left(s^2 \lambda^3 \int_{\Gamma_T} \xi \psi^2 \, \d S\, \d t + \lambda \int_{\Gamma_T} |\nabla_\Gamma \psi|^2 \, \d S\, \d t \right),\\
|\widetilde{B}_{3}| & \le C \left(s^2 \lambda^4 \int_{\gamma_T} \xi^2 \psi^2 \, \d S\, \d t + \lambda^2 \int_{\gamma_T} |\nabla_\Gamma \psi|^2 \, \d S\, \d t \right),
\end{align}
where we used the fact that $\nabla_\Gamma\eta=0$ on $\Gamma\setminus \gamma$ in the last estimate.
\bigskip

\noindent\textbf{Step 3. The transformed estimate.}
By summing the final inequalities obtained in previous steps, adding the complementary terms on $\gamma_T$, and choosing $\lambda_1, s_1$ large enough to absorb terms with lower powers in $s$ and $\lambda$, we arrive at
\begin{align*}
& \sum_{i,j} \langle (M_1\psi)_i, (M_2\psi)_j\rangle_{L^2(\Omega_T)}
    +  \sum_{i,j} \langle (N_1\psi)_i, (N_2\psi)_j\rangle_{L^2(\Gamma_T)}\\
& \quad  + C s^3\lambda^3 \int_{\gamma_T}  \xi^3 \psi^2  \, \d S\, \d t + C s\lambda \int_{\gamma_T} \xi(\partial_\nu \psi)^2  \, \d S\, \d t + C s\lambda \int_{\gamma_T}  \xi  |\nabla_\Gamma \psi|^2 \, \d S\, \d t\\
& \geq C s^3\lambda^4 \int_{\Omega_T} \xi^3 \psi^2 \, \d x\, \d t
+ C s^3\lambda^3 \int_{\Gamma_T}  \xi^3   \psi^2  \, \d S\, \d t 
  + C s\lambda^2 \int_{\Omega_T}  \xi   |\nabla \psi|^2  \, \d x\, \d t \\
  & \qquad + C s\lambda \int_{\Gamma_T} \xi(\partial_\nu \psi)^2  \, \d S\, \d t  + C s\lambda \int_{\Gamma_T}  \xi  |\nabla_\Gamma \psi|^2 \, \d S\, \d t\\
  &\qquad -Cs^3\lambda^4 \int_{\gamma_T}\xi^3\psi^2 \, \d S\, \d t -Cs\lambda^2 \int_{\gamma_T}\xi|\nabla_{\Gamma}\psi|^2 \, \d S\, \d t\\
 &\qquad - s^{-1} \int_{\Gamma_T} \xi^{-1} |\Delta_\Gamma \psi|^2 \, \d S\, \d t.
\end{align*}
The lower order terms coming from $\tilde f$ can be absorbed to the left-hand side by increasing $\lambda_1$ and $s_1$. The same thing holds true for $\tilde g$ except the boundary term $s^2\lambda^4 \displaystyle\int_{\gamma_T}\xi^2\psi^2 \, \d S\, \d t$, which remains on the right-hand side. Thus
\begin{align}
\|&M_1\psi\|_{L^2(\Omega_T)}^2 + \|M_2\psi\|_{L^2(\Omega_T)}^2  +\|N_1\psi\|_{L^2(\Gamma_T)}^2 
  + \|N_2\psi\|_{L^2(\Gamma_T)}^2 + s^3\lambda^4 \int_{\Omega_T} \xi^3 \psi^2 \, \d x\, \d t \notag\\
& +s\lambda^2 \!\int_{\Omega_T}\!  \xi   |\nabla \psi|^2  \, \d x\, \d t 
+  s^3\lambda^3\! \int_{\Gamma_T} \!\xi^3  \psi^2  \, \d S\, \d t 
  + s\lambda\! \int_{\Gamma_T}\! \xi  |\nabla_\Gamma\psi|^2  \, \d S\, \d t \notag\\
  & + s\lambda \!\int_{\Gamma_T}  \!  \xi(\partial_\nu \psi)^2  \, \d S\, \d t\notag\\
&\hspace{-0.2cm}\leq Cs^3 \lambda^4 \int_{\gamma_T} \xi^3\psi^2 \, \d x\, \d t
+ Cs\lambda^2 \int_{\gamma_T}\xi|\nabla_{\Gamma}\psi|^2 \, \d S\, \d t
 + C s\lambda \int_{\gamma_T} \xi(\partial_\nu \psi)^2  \, \d S\, \d t \notag\\
& \qquad  + s^{-1} \int_{\Gamma_T} \xi^{-1} |\Delta_\Gamma \psi|^2 \, \d S\, \d t .\label{m18}
\end{align}
To put the last summand in \eqref{m18} to the left-side, we make use of the equality
$\delta \Delta_\Gamma \psi = N_2\psi - s\psi \partial_t \alpha + \partial_\nu \psi-s^2\lambda^2 \xi^2\psi |\nabla_{\Gamma}\eta|^2$. Combined with 
\eqref{m6}, this equality yields
\begin{align}
I &:=s^{-1}\!\int_{\Gamma_T} \!\xi^{-1} |\Delta_\Gamma \psi|^2 \, \d S\, \d t\notag\\
&\leq \frac12\|N_2\psi\|_{L^2(\Gamma_T)}^2 + C s \! \int_{\Gamma_T}\! \xi^3 \psi^2 \, \d S\, \d t 
+ C\!\int_{\Gamma_T}\!\xi (\partial_\nu \psi)^2 \, \d S\, \d t \notag\\
& \quad +Cs^3\lambda^4\int_{\gamma_T}\xi^3 \psi^2 \, \d S\, \d t  \label{eq:LB}
\end{align}
for sufficiently large $s_1$.
We can now  choose sufficiently large $\lambda_1$ and $s_1$  
so that \eqref{m18} becomes
\begin{align}
\|&M_1\psi\|_{L^2(\Omega_T)}^2 + \|M_2\psi\|_{L^2(\Omega_T)}^2  +\|N_1\psi\|_{L^2(\Gamma_T)}^2 
 + \|N_2\psi\|_{L^2(\Gamma_T)}^2  + s^3\lambda^4 \int_{\Omega_T} \xi^3 \psi^2 \, \d x\, \d t \notag\\
 &  +s\lambda^2\!\int_{\Omega_T} \! \xi   |\nabla \psi|^2  \, \d x\, \d t
  +  s^3\lambda^3 \! \int_{\Gamma_T} \xi^3 \! \psi^2  \, \d S\, \d t 
   + s\lambda \! \int_{\Gamma_T}\!\xi(\partial_\nu \psi)^2 \, \d S\, \d t \notag\\
& + s\lambda \!\int_{\Gamma_T}\! \xi  |\nabla_\Gamma\psi|^2  \, \d S\, \d t\notag\\
& \leq Cs^3 \lambda^4 \int_{\gamma_T} \xi^3\psi^2 \, \d S\, \d t
+ Cs\lambda^2 \int_{\gamma_T}\xi|\nabla_{\Gamma}\psi|^2 \, \d S\, \d t + C s\lambda \int_{\gamma_T} \xi(\partial_\nu \psi)^2  \, \d S\, \d t. \label{est:car0}
\end{align}
\smallskip

\noindent\textbf{Step 4. Inverting the transformation.}
The inequality \eqref{eq:LB} allows to replace in 
\eqref{est:car0} the summand $\|N_2\psi\|_{L^2(\Gamma_T)}^2$ by the term $I$ times a constant, where we 
increase $\lambda_1$ and $s_1$ if necessary to absorb the lower order terms in \eqref{eq:LB}. Similarly,
from $$\partial_t\psi= N_1 \psi- s\lambda \psi \xi \partial_\nu \eta + 2s\lambda^2\xi \psi |\nabla_{\Gamma} \eta|^2 +2s\lambda \xi \nabla_{\Gamma}\psi \cdot\nabla_{\Gamma}\eta,$$ we deduce that 
\begin{align}
I' &:=s^{-1} \int_{\Gamma_T} \xi^{-1} |\partial_t\psi|^2 \, \d S\, \d t 
 \leq C\|N_1\psi\|_{L^2(\Gamma_T)}^2 + C s\lambda^2  \int_{\Gamma_T} \xi \psi^2 \, \d S\, \d t \notag\\
 &\qquad + Cs\lambda^4\int_{\gamma_T}\xi\psi^2 \, \d S\, \d t + Cs\lambda^2 \int_{\gamma_T}\xi|\nabla_{\Gamma}\psi|^2 \, \d S\, \d t,
 \end{align}
and then $I'$ can be inserted on the left-hand side of \eqref{est:car0}. Similarly, one can handle the corresponding terms on $\Omega_T$, (see \cite{FG'06}). We deduce that
\begin{align}
   &s^{-1}\int_{\Omega_T} \xi^{-1} (|\partial_t\psi|^2 + |\Delta \psi|^2 )\, \d x\, \d t  
  +s^{-1}\int_{\Gamma_T} \xi^{-1} (|\partial_t\psi|^2 + |\Delta_\Gamma \psi|^2) \, \d S\, \d t  \notag\\
 &\qquad  + s^3\lambda^4 \int_{\Omega_T} \xi^3 \psi^2 \, \d x\, \d t 
  +s\lambda^2 \int_{\Omega_T} \xi   |\nabla \psi|^2  \, \d x\, \d t \notag\\
 &\qquad  +  s^3\lambda^3 \int_{\Gamma_T} \xi^3  \psi^2  \, \d S\, \d t 
   + s\lambda \int_{\Gamma_T}\xi(\partial_\nu \psi)^2 \, \d S\, \d t
    + s\lambda \int_{\Gamma_T}\xi  |\nabla_\Gamma\psi|^2  \, \d S\, \d t\notag\\
& \leq Cs^3 \lambda^4 \int_{\gamma_T} \xi^3\psi^2 \, \d S\, \d t
+ Cs\lambda^2 \int_{\gamma_T}\xi|\nabla_{\Gamma}\psi|^2 \, \d S\, \d t + C s\lambda \int_{\gamma_T} \xi(\partial_\nu \psi)^2  \, \d S\, \d t.
     \label{est:car1} 
\end{align}
To estimate the penultimate integral in the right-hand side, we consider a cut-off function $\chi \in C^2(\Gamma^0)$
such that $\chi=1  \text { on } \gamma, 0 \leq \chi \leq 1 \text { on } \Gamma$. Then
\begin{align*}
s \lambda^{2} \int_{\gamma_T} \xi\left|\nabla_{\Gamma} \psi\right|^{2} \, \d S\, \d t & \leq s \lambda^{2} \int_{\Gamma_{T}^{0}} \chi \xi\left|\nabla_{\Gamma} \psi\right|^{2} \, \d S\, \d t \\
& =-s \lambda^{2} \int_{\Gamma_{T}^{0}} \chi \xi \psi \Delta_{\Gamma} \psi \, \d S\, \d t -s \lambda^{2} \int_{\Gamma_{T}^{0}} \xi \psi \nabla_{\Gamma} \chi \cdot \nabla_{\Gamma} \psi \, \d S\, \d t \\
&\quad -s \lambda^{3} \int_{\Gamma_{T}^{0}} \chi \xi \psi \nabla_{\Gamma} \eta \cdot \nabla_{\Gamma} \psi \, \d S\, \d t \\
&\leq \epsilon s^{-1} \int_{\Gamma_{T}^{0}} \xi^{-1}\left|\Delta_{\Gamma} \psi\right|^{2} \, \d S\, \d t + \dfrac{1}{4\epsilon} s^3 \lambda^4 \int_{\Gamma_{T}^{0}} \xi^{3}|\psi|^{2} \, \d S\, \d t \\
& \qquad + Cs \lambda^{4} \int_{\Gamma_{T}^{0}} \xi |\psi|^{2} \, \d S\, \d t.
\end{align*}
Choosing $\epsilon >0$ sufficiently small, we obtain that
\begin{align}
   &s^{-1}\int_{\Omega_T} \xi^{-1} (|\partial_t\psi|^2 + |\Delta \psi|^2 )\, \d x\, \d t  
  +s^{-1}\int_{\Gamma_T} \xi^{-1} (|\partial_t\psi|^2 + |\Delta_\Gamma \psi|^2) \, \d S\, \d t  \notag\\
 &\qquad  + s^3\lambda^4 \int_{\Omega_T} \xi^3 \psi^2 \, \d x\, \d t 
  +s\lambda^2 \int_{\Omega_T} \xi   |\nabla \psi|^2  \, \d x\, \d t \notag\\
 &\qquad  +  s^3\lambda^3 \int_{\Gamma_T} \xi^3  \psi^2  \, \d S\, \d t 
   + s\lambda \int_{\Gamma_T}\xi(\partial_\nu \psi)^2 \, \d S\, \d t
    + s\lambda \int_{\Gamma_T}\xi  |\nabla_\Gamma\psi|^2  \, \d S\, \d t\notag\\
& \leq Cs^3 \lambda^4 \int_{\gamma_T} \xi^3\psi^2 \, \d S\, \d t + C s\lambda \int_{\gamma_T} \xi(\partial_\nu \psi)^2  \, \d S\, \d t.
     \label{est:car2}
\end{align}
\smallskip

\noindent\textbf{Step 5. Estimating the normal derivative term.}
The most challenging term to be eliminated is the last term in \eqref{est:car2}. For this purpose, we proceed as follows: since
\begin{align*}
s\lambda \int_{\gamma_T}\xi(\partial_\nu \psi)^2 \,\d S\,\d t &= s\lambda \int_{\gamma_T} \xi\left(\mathrm{e}^{-s\alpha} \partial_\nu \varphi + s\lambda \xi \partial_\nu \eta \psi \right)^2 \,\d S\,\d t\\
& \le C s\lambda \int_{\gamma_T} \left(\mathrm{e}^{-2 s \alpha} \xi (\partial_\nu \varphi)^2 + s^2\lambda^2\xi^3 \psi^2 \right) \,\d S\,\d t\\
& \le C s\lambda \int_{\gamma_T} \mathrm{e}^{-2 s \alpha} \xi (\partial_\nu \varphi)^2 \,\d S\,\d t + C s^3\lambda^3 \int_{\gamma_T}\xi^3 \psi^2 \,\d S\,\d t,
\end{align*}
it suffices to estimate
$$
s \lambda \int_{\gamma_{T}} \mathrm{e}^{-2 s \alpha} \xi \left|\partial_{\nu} \varphi\right|^{2} \,\d S\,\d t.
$$
Since $d>0$, using $\eqref{eq:adj}_2$,  we have
\begin{equation}
s \lambda \int_{\gamma_{T}} \mathrm{e}^{-2 s \alpha} \xi\left|\partial_{\nu} \varphi\right|^{2} \,\d S\,\d t \leq C s \lambda \int_{\gamma_{T}} \mathrm{e}^{-2 s \alpha} \xi \left(\left|\Delta_{\Gamma} \varphi\right|^{2}+\left|\partial_{t} \varphi\right|^{2}\right) \,\d S\,\d t. \label{eQ1}
\end{equation}
Next, we estimate the two terms in the right-hand side of \eqref{eQ1}. We introduce a cut-off function $\chi \in C^2(\Gamma^0)$
such that $0 \leq \chi \leq 1 \text { on } \Gamma$ and $\chi=1 \text { on } \gamma$. Then, using the surface divergence formula twice, we obtain
\begin{align}
& s \lambda \int_{\gamma_{T}} \mathrm{e}^{-2 s \alpha} \xi \left|\Delta_{\Gamma} \varphi\right|^{2} \,\d S\,\d t \le s \lambda \int_{\Gamma_{T}^0} \mathrm{e}^{-2 s \alpha} \chi \xi \left|\Delta_{\Gamma} \varphi\right|^{2} \,\d S\,\d t \notag\\
& \le s\lambda \int_{\Gamma_{T}^0} \varphi \Delta_\Gamma\left(\mathrm{e}^{-2 s \alpha} \chi \xi \Delta_\Gamma \varphi\right) \,\d S\,\d t \notag \\
&  = \int_{\Gamma_{T}^0} \mathrm{e}^{-2 s \alpha} \varphi \Big(\left(6 s^2\lambda^3 \xi^2 + 4 s^3\lambda^3 \xi^3 + s\lambda^3 \xi\right) (\chi |\nabla_\Gamma \eta|^2 \Delta_\Gamma \varphi) \notag \\
& \qquad + (2 s^2\lambda^2 \xi^2 + s\lambda^2 \xi)(\chi \Delta_\Gamma \eta \Delta_\Gamma \varphi)
+ s\lambda \xi \Delta_\Gamma \chi \Delta_\Gamma \varphi + s\lambda \xi \chi \Delta_\Gamma^2 \varphi\notag \\
& \qquad + (4 s^2\lambda^2 \xi^2 + 2 s\lambda^2\xi) (\nabla_\Gamma \eta \cdot \nabla_\Gamma \chi) \Delta_\Gamma \varphi \notag\\
&\qquad + (4 s^2\lambda^2 \xi^2 + 2 s\lambda^2\xi) \chi\nabla_\Gamma \eta \cdot \nabla_\Gamma \Delta_\Gamma \varphi + 2 s\lambda \xi \nabla_\Gamma \chi \cdot \nabla_\Gamma \Delta_\Gamma \varphi \Big) \,\d S\,\d t \notag\\
& \le C \left(s^{6} \lambda^{6} \int_{\Gamma_{T}^{0}} \mathrm{e}^{-(2-3 \mu) s \alpha} \xi^{6} |\varphi|^2 \,\d S\,\d t +\int_{0}^{T} \max _{x \in \overline{\Omega}}\left(\mathrm{e}^{-(2+3 \mu) s \alpha}\right) \|\varphi\|_{H^{4}(\Gamma)}^{2} \,\d t \right) \notag
\end{align}
for all $\mu>0$, where we used Young's inequality.
Integrating by parts in time, and using $\xi \le CT^2 \xi^2, \,|\xi_t| + |\alpha_t| \le C T\xi^2, \, |\xi_{tt}|+|\alpha_{tt}| \le C T^2 \xi^3$ together with the fact that $\mathrm{e}^{-2 s \alpha}\rvert_{t=0}=\mathrm{e}^{-2 s \alpha}\rvert_{t=T}=0$, we obtain
\begin{align}
& s\lambda \int_{\Gamma_{T}^0} \mathrm{e}^{-2 s \alpha} \xi |\partial_t \varphi|^2 \,\d S\,\d t =  - s\lambda \int_{\Gamma_{T}^0} \varphi \partial_t \left(\mathrm{e}^{-2 s \alpha} \xi \partial_t \varphi\right) \,\d S\,\d t \notag\\
& = s \lambda \int_{\Gamma_{T}^0} \mathrm{e}^{-2 s \alpha}\left(s\alpha_t \xi - \frac{1}{2} \xi_t\right)  \partial_t (\varphi^2) \,\d S\,\d t - s\lambda \int_{\Gamma_{T}^0} \mathrm{e}^{-2 s \alpha} \xi \varphi \varphi_{tt} \,\d S\,\d t \notag\\
& = s \lambda \int_{\Gamma_{T}^0} \mathrm{e}^{-2 s \alpha} \left(2 s^2\alpha_t^2 \xi -2 s\alpha_t \xi_t - s\alpha_{tt} \xi + \frac{1}{2} \xi_{tt}\right)  \varphi^2 \,\d S\,\d t \notag\\
& \quad - s\lambda \int_{\Gamma_{T}^0} \mathrm{e}^{-2 s \alpha} \xi \varphi \varphi_{tt} \,\d S\,\d t \notag\\
& \le C \left(s^{6} \lambda^{6} \int_{\Gamma_{T}^{0}} \mathrm{e}^{-(2-3 \mu) s \alpha} \xi^{6} |\varphi|^2 \,\d S\,\d t +\int_{0}^{T} \max _{x \in \overline{\Omega}}\left(\mathrm{e}^{-(2+3 \mu) s \alpha}\right) \|\varphi_{tt}\|_{L^2(\Gamma)}^{2} \,\d t \right).
\end{align}
Then
\begin{align}
& s \lambda \int_{\gamma_{T}} \mathrm{e}^{-2 s \alpha} \xi\left|\partial_{\nu} \varphi\right|^{2} \,\d S\,\d t \le C\left(s^{6} \lambda^{6} \int_{\Gamma_{T}^{0}} \mathrm{e}^{-(2-3 \mu) s \alpha} \xi^{6} \left|\varphi_{\Gamma}\right|^{2} \,\d S\,\d t \right. \notag\\
&\hspace{4.2cm}\left.  +\int_{0}^{T} \max _{x \in \overline{\Omega}}\left(\mathrm{e}^{-(2+3 \mu) s \alpha}\right)\left(\|\varphi\|_{H^{4}(\Gamma)}^{2}+\left\|\varphi_{t t}\right\|_{L^{2}(\Gamma)}^{2}\right) \,\d t\right) .
\end{align}
Moreover, taking $\lambda \geq C_{\mu}$ large enough, we obtain $$\max\limits_{x \in \overline{\Omega}}\left(\mathrm{e}^{-(2+3 \mu) s \alpha}\right) \le \min\limits_{x \in \overline{\Omega}}\left(\mathrm{e}^{-(2+2 \mu) s \alpha}\right).$$
Let us denote $\alpha^{*}(t)=\max\limits_{x \in \overline{\Omega}} \alpha(t, x)$. Then $\min\limits_{x \in \overline{\Omega}}\left(\mathrm{e}^{-(2+2 \mu) s \alpha}\right)=\mathrm{e}^{-(2+2 \mu) s \alpha^*}$, and it suffices to estimate the quantity
$$
\int_{0}^{T} \left(\|q\|_{H^{4}(\Gamma)}^{2}+\left\|q_{t t}\right\|_{L^{2}(\Gamma)}^{2} \right) \,\d t,
$$
where $(q, q_\Gamma)=\mathrm{e}^{-(1+\mu) s \alpha^*} (\varphi, \varphi_\Gamma)$, by using parabolic regularity estimates. Remark that $(q, q_\Gamma)$ satisfies the system
\begin{empheq}[left = \empheqlbrace]{alignat=2}
\begin{aligned}
&-\partial_t q -d\Delta q  = (1+\mu) s \alpha_t^* q, &\text{ in } \Omega_T , \\
&-\partial_t q_\Gamma -\delta\Delta_\Gamma q_\Gamma +d\partial_\nu q  = (1+\mu) s \alpha_t^* q_\Gamma, &\qquad\text{ on } \Gamma_T , \\
& q_{\Gamma}(t,x) = q_{|\Gamma}(t,x), &\qquad\text{ on } \Gamma_T , \\
& (q(T,\cdot),q_\Gamma(T,\cdot)) = (0, 0), &\Omega\times\Gamma. \label{eq:q}
\end{aligned}
\end{empheq}
We use Proposition \ref{prop:accreg}, $(iii)$ for \eqref{eq:q}, and then $(ii)$ for the system satisfied by $\alpha_t^*(q,q_\Gamma)$, we obtain
\begin{align}
& \hspace{0.5cm}\int_{0}^{T} \left(\|q\|_{H^{4}(\Gamma)}^{2}+\left\|q_{t t}\right\|_{L^{2}(\Gamma)}^{2} \right) \,\d t \notag\\
& \le C_{\mu}(1+T)^2 \left(\|s \alpha_t^* (q,q_\Gamma)\|^2_{L^2(0,T;\mathbb{H}^2)} +  \|s \partial_t \left(\alpha_t^* (q,q_\Gamma)\right)\|^2_{L^2(0,T;\mathbb{L}^2)}\right)\notag\\
&\le C_{\mu}(1+T)^4 \|s \left(\alpha_{tt}^*-s(1+\mu)(\alpha_t^*)^2\right) (q,q_\Gamma)\|^2_{L^2(0,T; \mathbb{L}^2)} \notag\\
& \le C_{\mu} T^4(1+T)^4 \left(\int_{\Omega_T} \mathrm{e}^{-(2+\mu) s \alpha^*} |\varphi|^2 \,\d x\,\d t + \int_{\Gamma_T} \mathrm{e}^{-(2+\mu) s \alpha^*} |\varphi_\Gamma|^2 \,\d S\,\d t\right),
\end{align}
where we used $(s\alpha_{tt}^*)^2+(s\alpha_t^*)^4 \le C T^4 (s\alpha^*)^8 \le C T^4 \mathrm{e}^{\mu s \alpha^*}$. Taking $\mu=\frac{1}{3}$, the previous integrals can be absorbed by the left-hand side of the inequality \eqref{est:car2}, and the rest of the proof follows from the same strategy as in \cite{FG'06}.
\end{proof}

\section{Boundary null controllability}\label{sec:null}
In this section, we prove the null controllability of the following linear inhomogeneous system
\begin{empheq}[left = \empheqlbrace]{alignat=2}
\begin{aligned}
&\partial_t y - d\Delta y = f(t,x), &\qquad\text{ in } \Omega_T , \\
&\partial_t y_{\Gamma} -\delta\Delta_\Gamma y_\Gamma + d\partial_\nu y = v(t,x) \mathds{1}_{\Gamma^0} + g(t,x), &\qquad\text{ on } \Gamma_T , \\
& y_{\Gamma}(t,x) = y_{|\Gamma}(t,x), &\qquad\text{ on } \Gamma_T , \\
& (y,y_\Gamma)\rvert_{t=0}=(y_0,y_{0,\Gamma}), &\Omega\times\Gamma. \label{eq:ncont}
\end{aligned}
\end{empheq}
To this end, we shall first infer an observability estimate for the homogeneous backward system given by
\begin{empheq}[left = \empheqlbrace]{alignat=2}
\begin{aligned}
&-\partial_t \varphi-d\Delta \varphi = 0,  &\qquad \text{in }\Omega_T, \\
&-\partial_t \varphi_\Gamma-\delta\Delta_\Gamma \varphi_\Gamma+d\partial_\nu \varphi =0,  &\qquad \text{on }\Gamma_T, \\
&\varphi_{|\Gamma}(t,x)=\varphi_{|\Gamma}(t,x), &\qquad\text{on } \Gamma_T, \\
&(\varphi(T,\cdot),\varphi_\Gamma(T,\cdot)) = (\varphi_T, \varphi_{T,\Gamma}), &\Omega\times \Gamma, 
\label{eq:nadj}
\end{aligned}
\end{empheq}
by applying the Carleman estimate \eqref{Car}.
\begin{proposition}\label{lem:observa}
There exists a constant $C=C\left(\Omega,\Gamma^0\right) >0$ such that,  for all $(\varphi_T, \varphi_{T,\Gamma}) \in \mathbb{L}^2$, the mild solution 
$(\varphi, \varphi_{\Gamma})$ of the backward problem \eqref{eq:nadj} satisfies
\begin{equation*}
 \|\varphi(0,\cdot)\|^2_{L^2(\Omega)} +\|\varphi_\Gamma(0,\cdot)\|^2_{L^2(\Gamma)} \leq \mathrm{e}^{C\left(T^{2/3} +\frac{1}{T}\right)}\int_{\Gamma^0_T} |\varphi|^2 \, \d S\, \d t.
\end{equation*}
Moreover, for $ (y_0, y_{0,\Gamma}) \in \mathbb{L}^2$ the mild solution $(y,y_\Gamma)$ of the forward problem \eqref{eq:ncont} with $f=g=v=0$
satisfies
\begin{equation*}
 \|y(T,\cdot)\|^2_{L^2(\Omega)} +\|y_\Gamma(T,\cdot)\|^2_{L^2(\Gamma)} \leq \mathrm{e}^{C\left(T^{2/3} +\frac{1}{T}\right)} \int_{\Gamma^0_T} |y|^2 \, \d S\, \d t. 
\end{equation*}
\end{proposition}

\begin{proof} 
The proof mainly follows from the Carleman estimate \eqref{Car} with an adaptation of the same ideas in \cite[Proposition 5]{KM'19}, so we omit the details.
\end{proof}

\begin{remark}
We emphasize that the boundary observability constant
$$K_{\Gamma^0}\left(\Omega,T\right):=\exp\left({C\left(\Omega,\Gamma^0\right)\left(T^{2/3} +\frac{1}{T}\right)}\right)$$
obtained in Proposition \ref{lem:observa} is different from the interior observability constant
$$K_\omega\left(\Omega, T\right):=\exp\left({C\left(\Omega,\omega\right)\left(1+\frac{1}{T}\right)}\right)$$
obtained in \cite[Proposition 5]{KM'19}. This is merely due to the parabolic regularity estimates used in the Carleman estimate to absorb the normal derivative term.
\end{remark}

Now, we are in position to conclude the null controllability of \eqref{eq:ncont}, where inhomogeneities with exponential decay at $t=0$ and $t=T$ are allowed. Recall that
$$
\theta(t) =\frac{1}{t(T-t)},  \quad\qquad t\in (0,T).
$$
Let us introduce the $L^2$-weighted spaces
\begin{align*}
Z_{\Omega} &:= \left\{f\in L^2(\Omega_T): \mathrm{e}^{s\theta(t)} \theta^{-\frac{3}{2}}(t)\, f\in L^2(\Omega_T)\right\},\\
Z_{\Gamma} &:= \left\{g\in L^2(\Gamma_T): \mathrm{e}^{s\theta(t)} \theta^{-\frac{3}{2}}(t)\, g\in L^2(\Gamma_T)\right\},\\
X &:= \left\{Y\in L^2\left(0,T;\mathbb{L}^2\right): \mathrm{e}^{s t\theta(t)} Y\in L^2\left(0,T;\mathbb{L}^2\right)\right\},
\end{align*}
with the inner products
\begin{align*}
\langle f_1,f_2\rangle_{Z_{\Omega}} &= \int_{\Omega_T} f_1f_2 \,\mathrm{e}^{2s\theta} \, \theta^{-3}\, \d x\, \d t \\
\langle g_1,g_2\rangle_{Z_{\Gamma}} &= \int_{\Gamma_T} g_1g_2 \,\mathrm{e}^{2s\theta}\,\theta^{-3} \, \d S\, \d t.
\end{align*}

\begin{theorem}\label{thm:control-linear} 
For any initial data $Y_0:=(y_0, y_{0,\Gamma})\in \mathbb{L}^2$, $f\in Z_\Omega$ and $g\in Z_\Gamma,$ there exists a control $v\in L^2\left(\Gamma^0_T\right)$ such that the mild solution $(y,y_\Gamma)$ of 
\eqref{eq:ncont} satisfies $$y(T,\cdot) = y_\Gamma(T,\cdot)=0.$$ Moreover, there exists a constant $C>0$ independent of $Y_0$, $f$ and $g$ such that
\begin{equation}\label{contest}
\|Y\|_X + \|v\|_{L^2\left(\Gamma^0_T\right)} \le C \left(\|Y_0\|_{\mathbb{L}^2}+ \|f\|_{Z_\Omega}+\|g\|_{Z_\Gamma}\right).
\end{equation}
\end{theorem}

\begin{proof} 
By duality, the proposition will be derived from the previous observability inequality along with the Carleman estimate. To do this, we define a bounded linear operator
$\mathcal T: L^2\left(\Gamma^0_T\right)\to \mathbb L^2$ by 
\[ \mathcal Tv= \int_0^T \mathrm{e}^{(T-\tau)\mathcal{A}}(0,-\mathds{1}_{\Gamma^0} v(\tau))\, \d\tau.\]
Using the embedding
$Z_\Omega\times Z_\Gamma\hookrightarrow L^2(\Omega_T) \times L^2(\Gamma_T)$, we consider
the bounded linear operator $\mathcal S: \mathbb L^2\times Z_\Omega\times Z_\Gamma \to \mathbb{L}^2$ defined as follows
\[\mathcal S(Y_0, f, g) = \mathrm{e}^{T\mathcal{A}}Y_0+ \int_0^T \mathrm{e}^{(T-\tau)\mathcal{A}}(f(\tau),g(\tau))\, \d\tau.  \]
By Proposition~\ref{prop:mild}\;(a), the function $\mathcal S (Y_0,f,g) - \mathcal T v $ is the final value $(y(T),y_\Gamma(T))$ of the mild solution $(y,y_\Gamma)$ to the inhomogeneous system \eqref{eq:ncont}. Note that the adjoint operator $\mathcal T^*: \mathbb L^2 \to L^2\left(\Gamma^0_T\right)$ is given by 
\[ \mathcal T^* (\varphi_T, \varphi_{T,\Gamma}) = -\mathds{1}_{\Gamma^0} \varphi_{\Gamma},\] 
where $(\varphi(t), \varphi_{\Gamma}(t))= \mathrm{e}^{(T-t)\mathcal{A}}(\varphi_T, \varphi_{T,\Gamma})$ is the mild solution of the homogeneous backward problem \eqref{eq:nadj} with final value $(\varphi_T, \varphi_{T,\Gamma})$, see Proposition~\ref{prop:mild}\,(d).\\
The adjoint operator $\mathcal S^*: \mathbb L^2 \to 
 \mathbb{L}^2\times Z_\Omega\times Z_\Gamma$ is given by
\[\mathcal S^* (\varphi_T, \varphi_{T,\Gamma})= \left((\varphi(0), \varphi_{\Gamma}(0)),\mathrm{e}^{-2s\theta} \theta^3 \varphi, \mathrm{e}^{-2s\theta} \theta^3 \varphi_{|\Gamma}\right).\]
From Proposition~\ref{prop:mild}, we see that $(\varphi, \varphi_{\Gamma})\in\mathbb{E}_1$ is a strong solution of the backward
problem if $(\varphi_T, \varphi_{T,\Gamma})\in \mathbb{H}^1$. The observability inequality of Proposition~\ref{lem:observa} and the Carleman estimate from Lemma~\ref{lem:carleman} imply that
\begin{align*}
\|\mathcal S^* (\varphi_T, \varphi_{T,\Gamma})
\|_{\mathbb L^2\times Z_\Omega\times Z_\Gamma}^2 
& \le \|\varphi(0,\cdot)\|^2_{L^2(\Omega)} +\|\varphi_\Gamma(0,\cdot)\|^2_{L^2(\Gamma)}  + C\int_{\Omega_T} \mathrm{e}^{-2s\alpha} \xi^3 \varphi^2 \, \d x\, \d t \\
  & \qquad  + C\int_{\Gamma_T} \mathrm{e}^{-2s\alpha} \xi^3 \varphi^2 \, \d S\, \d t\\
& \leq C \int_{\Gamma^0_T} |\varphi|^2 \, \d S\, \d t = C \|\mathcal T^* (\varphi_T, \varphi_{T,\Gamma})\|_{L^2\left(\Gamma^0_T\right)}
\end{align*} 
for $(\varphi_T, \varphi_{T,\Gamma})\in \mathbb{H}^1$, and hence for $(\varphi_T, \varphi_{T,\Gamma})\in \mathbb{L}^2$ by denseness.  
Theorem~15.2 of \cite{Za'20} implies that the range of $\mathcal T$ contains that of $\mathcal S$; i.e.,
for all $Y_0\in \mathbb L^2$, $f\in Z_\Omega$ and $g\in Z_\Gamma$ there exists a control $v\in L^2\left(\Gamma^0_T\right)$ such that $\mathcal S (Y_0,f,g) = \mathcal T v$. Thus,  
$(y(T),y_\Gamma(T)) = \mathcal S (Y_0,f,g) - \mathcal T v = (0,0)$. Finally, the estimate \eqref{contest} follows by arguing as in \cite[Proposition 4.3]{MMS'17}
\end{proof}

\begin{remark}
Since the operator $\mathcal{A}$ is self-adjoint dissipative and the system \eqref{eq:ncont} without source terms is null controllable, a similar controllability result in the presence of source terms can be obtained from \cite[Proposition 2.3]{LTT'13}. In the latter result, the source terms $f$ and $g$ should be taken respectively in the weighted spaces $\tilde{Z}_{\Omega}$ and $\tilde{Z}_{\Gamma}$ which result by considering a weight function $(\rho(t))^{-1}$ instead of $\mathrm{e}^{s\theta(t)} \theta^{-\frac{3}{2}}(t)$ in $Z_{\Omega}$ and $Z_{\Gamma}$, where $\rho$ is an appropriately chosen function.
\end{remark}

\section*{Acknowledgment} The authors would like to thank Jon Asier Bárcena-Petisco for fruitful discussions and invaluable comments.

\end{document}